\documentclass[12pt]{amsart} \textwidth=14.5cm \oddsidemargin=1cm
\usepackage{amsmath}
\usepackage{color}
\usepackage{amsxtra}
\usepackage{amscd}
\usepackage{amsthm}
\usepackage{amsfonts}
\usepackage{amssymb}
\usepackage{eucal}

\input prepictex
\input pictex
\input postpictex

\newtheorem{thm}{Theorem}[section]
\newtheorem{lem}{Lemma}[section]

\newtheorem{prop}{Proposition}[section]

\theoremstyle{definition}
\newtheorem{definition}{Definition}[section]

\theoremstyle{remark}
\newtheorem{rem}{Remark}[section]

\numberwithin{equation}{section}

\begin{document}

\newcommand{\thmref}[1]{Theorem~\ref{#1}}
\newcommand{\secref}[1]{Section~\ref{#1}}
\newcommand{\lemref}[1]{Lemma~\ref{#1}}
\newcommand{\propref}[1]{Proposition~\ref{#1}}
\newcommand{\corref}[1]{Corollary~\ref{#1}}
\newcommand{\remref}[1]{Remark~\ref{#1}}
\newcommand{\eqnref}[1]{(\ref{#1})}
\newcommand{\exref}[1]{Example~\ref{#1}}

\newcommand{\nc}{\newcommand}
\nc{\Z}{{\mathbb Z}}
\nc{\C}{{\mathbb C}}
\nc{\N}{{\mathbb N}}
\nc{\F}{{\mf F}}
\nc{\Q}{\ol{Q}}
\nc{\la}{\lambda}
\nc{\ep}{\epsilon}
\nc{\h}{\mathfrak h}
\nc{\n}{\mf n}
\newcommand{\glpm}{ \widehat{{ \mathfrak {gl} }}_{\pm} }
\newcommand{\hL}{{\Lambda} }

\nc{\G}{{\mathfrak g}}
\nc{\SG}{\overline{\mathfrak g}}
\nc{\Li}{{\mc L}} \nc{\La}{\Lambda} \nc{\is}{{\mathbf i}}
\nc{\V}{\mf V} \nc{\bi}{\bibitem} \nc{\NS}{\mf N}
\nc{\dt}{\mathord{\hbox{${\frac{d}{d t}}$}}} \nc{\E}{\mc E}
\nc{\ba}{\tilde{\pa}} \nc{\half}{\frac{1}{2}} \nc{\mc}{\mathcal}
\nc{\mf}{\mathfrak} \nc{\ol}{\fracline} \nc{\hf}{\frac{1}{2}}
\nc{\eps}{\epsilon}
\nc{\ohgl}{\widehat{\mathfrak{gl}}}
\nc{\hgl}{a_\infty}
\nc{\gl}{{\mathfrak{gl}}}
\nc{\hz}{\hf+\Z}
\nc{\vac}{|0 \rangle}
\nc{\dinfty}{{\infty\vert\infty}}
\nc{\SLa}{\overline{\Lambda}}
\nc{\SF}{\overline{\mathfrak F}}
\nc{\SP}{\overline{\mathcal P}}

\newcommand{\A}{\mathsf A}
\newcommand{\D}{\mathsf D}
\newcommand{\Ct}{\mathsf C}
\newcommand{\no}{\text{:}}
\newcommand{\trace}{\,{\rm tr}\,}
\newcommand{\ainf}{a_\infty}
\newcommand{\cinf}{c_\infty}
\newcommand{\dinf}{d_\infty}

\advance\headheight by 2pt

\title[Bloch-Okounkov functions of negative levels]
{The Bloch-Okounkov correlation functions of negative levels}

\author[Cheng]{Shun-Jen Cheng}
\address{Institute of Mathematics, Academia Sinica, Taipei, Taiwan 11529}
\email{chengsj@math.sinica.edu.tw}

\author[Taylor]{David G. Taylor}
\address{Department of Mathematics, University of Virginia, Charlottesville, VA
22904;
{\em (New address after August 2007)} Department of Mathematics,
Computer Science, and Physics, Roanoke College,
Salem, VA 24153} \email{taylor@roanoke.edu}

\author[Wang]{Weiqiang Wang}
\address{Department of Mathematics, University of Virginia, Charlottesville, VA 22904}
 \email{ww9c@virginia.edu}

\begin{abstract}
Bloch and Okounkov introduced an $n$-point correlation function on
the fermionic Fock space and found a closed formula in terms of
theta functions. This function affords several distinguished
interpretations and in particular can be formulated as correlation
functions on irreducible $\ohgl_\infty$-modules of level one. These
correlation functions have been generalized for irreducible
integrable modules of $\ohgl_\infty$ and its classical Lie
subalgebras of positive levels by the authors. In this paper we
extend further these results and compute the correlation functions
as well as the $q$-dimensions for modules of $\ohgl_\infty$ and its
classical subalgebras at negative levels.
\end{abstract}

\maketitle

\section{Introduction}

\subsection{}
Bloch-Okounkov \cite{BO} (also \cite{Ok}) formulated an $n$-point
correlation function on the fermionic Fock space and found a
beautiful closed formula for it in terms of Jacobi theta functions.
This function has since made appearances in several distinct setups
including Gromov-Witten theory and Hilbert schemes. The viewpoint
taken in that paper is to regard the original Bloch-Okounkov
functions as correlation functions on irreducible
$\ohgl_\infty$-modules of level one. In this sense, these
correlation functions have been generalized for irreducible
integrable modules of $\ohgl_\infty$ at positive levels \cite{CW},
and also generalized for integrable modules of classical Lie
subalgebras of $\ohgl_\infty$ (of type $B,C,D$) at positive levels
\cite{TW, W2} (also cf. \cite{Mil}). These Bloch-Okounkov functions
can be also viewed as a refined version of character formulas for
the corresponding modules.

The Lie algebra $\ohgl_\infty$ and its classical subalgebras were
introduced by the Kyoto school \cite{DJKM1, DJKM2} in connection
with vertex operators and KP integrable hierarchies, etc., and they
have played fundamental roles in representation theory of
infinite-dimensional Lie algebras. They are also intimately related
to the $W_{1+\infty}$ algebra and its classical subalgebras (cf.
\cite{KWY} and references therein).

\subsection{}
The goal of this paper is to study the Bloch-Okounkov correlation
functions for (mostly) irreducible highest weight modules of
negative levels. There are several dualities \cite{W1} (also see
\cite{KR} for type $A$) on various Fock spaces of bosonic ghosts
(also called $\beta\gamma$ systems in physics literature) between a
finite-dimensional Lie group on one side and $\ohgl_\infty$, or one
of its classical subalgebras, on the other side. These dualities are
natural infinite-dimensional generalizations of the classical Howe
duality for finite-dimensional Lie groups/algebras \cite{H1, H2}.
The modules of $\ohgl_\infty$ and other Lie algebras under
consideration in the present paper appear in these Howe duality
decompositions. In fact, these dualities are the essential tools
that allow us to reduce the calculation of the correlation functions
of a general negative level to those of the bottom levels (i.e.
level $-1$ or $-\hf$), and the precise relation involves summation
over the Weyl group of the corresponding Lie group in Howe duality.
Similar ideas have been used in \cite{CW, TW}.

We show that the $n$-point correlation functions at the bottom
levels satisfy certain $q$-difference equations (such an idea goes
back to \cite{BO, Ok} in the original setup). We are able to compute
the $1$-point and $2$-point correlation functions at the bottom
levels explicitly in terms of certain $q$-hypergeometric series (cf.
\cite{GR}), but the general $n$-point case remains open. The
negative level case is technically more complicated than the
positive level case treated in our earlier works \cite{CW, TW}, and
the difference between negative level and positive level is already
apparent at their bottom levels.

The modules considered in this paper all possess a natural
$\Z_+$-grading with finite-dimensional subspaces, and thus it
makes sense to ask for their $q$-dimension (i.e. graded
dimension). The strategy used to calculate the correlation
functions allows us to determine explicitly the $q$-dimensions
(which can be regarded as the $0$-point correlation function) for
the corresponding modules. The $q$-dimension formula at the bottom
level has been folklore, but the general case appears to be new.
Indeed, the $q$-dimension formula at the bottom level boils down
to an intriguing $q$-series identity, which affords several
different proofs to date \cite{BCMN, FeF, K1}. Each of these
proofs is complicated yet very interesting in its own way, using
{\em super} boson-fermion correspondence or the underlying
Virasoro algebra structure of the Fock space of bosonic ghosts,
just as the celebrated Jacobi triple product identity underlies
the boson-fermion correspondence. Here we offer a very short
combinatorial proof of this remarkable identity.

The consideration of the $\cinf$-modules of level $l-\hf$ in this
paper, which is the only positive level case left out of \cite{TW}
since it involves bosonic Fock space, also helps to complete the
study of the correlation functions and $q$-dimensions for integrable
modules of all classical Lie subalgebras of $\ohgl_\infty$.

\subsection{}
This paper is organized as follows.  We treat the $n$-point
correlation functions for $\ohgl_\infty$-modules of negative level
in Section~\ref{sec:ainf}. The same is then done in
Section~\ref{sec:cinf} for the classical Lie subalgebra of
$\ohgl_\infty$ of type $C$ and in Section~\ref{sec:dinf} for type
$D$. We note that the subalgebra of type $B$ of negative level does
not feature in a Howe duality and thus does not appear in this work.
Along the way, $q$-dimension formulas are also given in each case.

\subsection{Notations} For a classical simple Lie algebra, we use the
standard notation to denote the roots by
$\varepsilon_i-\varepsilon_j$, $\pm\varepsilon_i$ and $\pm
2\varepsilon_i$ etc. By $(\cdot,\cdot)$ we mean the usual
symmetric bilinear form determined by
$(\varepsilon_i,\varepsilon_j)=\delta_{ij}$ and write
$\|x\|^2=(x,x)$. Further we let
$\rho=\sum_{i=1}^l(l-i)\varepsilon_i$,
$\rho_B=\sum_{i=1}^l(l-i+\hf)\varepsilon_i$, and
$\rho_C=\sum_{i=1}^l(l-i+1)\varepsilon_i$.

Given a Lie algebra $x_\infty$ with $x=a,c,d$, we denote
$L(x_\infty; \Lambda,k)$ the irreducible $x_\infty$-module of
highest weight $\Lambda$ and level $k$.

We denote by $\N$ the set of natural numbers and by $\Z_+$ the set
of nonnegative integers.

{\bf Acknowledgment.} S-J.C. is partially supported by an NSC and an
Academia Sinica Investigator grant. He also thanks the Department of
Mathematics of the University of Virginia (UVa) for hospitality and
support. D.T. is supported by a Dissertation Semester Fellowship
from the Department of Mathematics at UVa. W.W. is partially
supported by NSF and NSA grants. We are thankful to the very capable
referee for corrections and suggestions for improvements.

\section{The $\ainf$-correlation functions and $q$-dimension formulas}
\label{sec:ainf}

\subsection{Some $q$-series identities}
We start with some combinatorial preparation. For an indeterminate
$q$ we let
$$\begin{aligned}
(a)_0&=1,\\
(a)_n&=(1-a)(1-aq)\cdots(1-aq^{n-1}),  \quad \text{ for } n\in\N, \\
(a)_\infty &= (1-a)(1-aq)(1-aq^2)\cdots = \prod_{i=0}^\infty
(1-aq^i).
\end{aligned}
$$
Alternatively we may also regard $q$ as a complex number with
$|q|<1$ to ensure the functions in this paper converge as analytic
functions. For $r,s\in\Z_+$, recall the $q$-hypergeometric series
(cf.~\cite{GR})
\begin{equation*}
_{r}\Phi_s(a_1,\cdots,a_{r};b_1,\cdots,b_s;z)=\sum_{n=0}^\infty
\frac{(a_1)_n\cdots(a_{r})_n}{(b_1)_n\cdots(b_s)_n(q)_n}\left((-1)^nq^{n(n-1)/2}\right)^{1+s-r}z^n.
\end{equation*}

In the sequel we will frequently use the following well-known
identities (cf. \cite{GR})
\begin{align}
  \sum_{m =0}^\infty \frac{(-z)^m q^{m(m-1)/2}}{(q)_m}
  = (z)_\infty,
  \qquad
  \sum_{l = 0}^\infty \frac{  (a)_l}{(q)_l} z^l
  = \frac{(az)_\infty}{(z)_\infty}.\label{exponential}
\end{align}

\begin{prop}\label{111}
For a given $k \ge 0$, we have
\begin{eqnarray}
\sum_{l=0}^\infty \frac{q^l}{(q)_{l}(tq)_{l+k}}
 = \frac1{(q)_\infty (tq)_\infty}
 \sum_{m \ge 0} (-1)^m q^{m(m+1)/2 +km} t^m.
\end{eqnarray}
\end{prop}

\begin{proof}
We calculate using (\ref{exponential})
{\allowdisplaybreaks
\begin{eqnarray*}
 (tq)_\infty \sum_{l=0}^\infty \frac{q^l}{(q)_{l}(tq)_{l+k}}
 &=& \sum_{l=0}^\infty \frac{q^l (tq^{l+1+k})_{\infty}}{(q)_{l} } \\
 &=& \sum_{l=0}^\infty \frac{q^l}{(q)_{l} }
   \sum_{m =0}^\infty \frac{(-tq^{l+1+k})^m q^{m(m-1)/2}}{(q)_m}\\
  &=& \sum_{m =0}^\infty \frac{(-t)^m q^{m(m+1)/2+mk}}{(q)_m}
    \sum_{l=0}^\infty \frac{q^{(m+1)l}}{(q)_{l} }\\
  &=& \sum_{m =0}^\infty \frac{(-t)^m q^{m(m+1)/2+mk}}{(q)_m}
  \frac{1}{(q^{m+1})_{\infty} }\\
  &=& \frac1{(q)_\infty}
 \sum_{m \ge 0} (-1)^m q^{m(m+1)/2+mk} t^m.
\end{eqnarray*}}
This finishes the proof.
\end{proof}

We give a short elementary proof of the following identity that
appeared in the literature with different but more complicated
proofs (cf. e.g. \cite{BCMN, FeF, K1}).

\begin{thm} \label{identity:ff}
We have
\begin{eqnarray*}
 \frac1{(u)_{\infty}(u^{-1}q)_{\infty}}
 &=& \frac1{(q)_{\infty}^2}   \sum_{m=0}^\infty (-1)^m q^{\hf m(m+1)}
   \left( \sum_{k\ge 0}   q^{km} u^k
   + \sum_{k >0} q^{k(m+1)} u^{-k} \right)  \\
    &=&  \frac1{(q)_{\infty}^2} \sum_{m \in \Z} (-1)^m q^{\hf m(m+1)} \frac1{1-uq^m}.
\end{eqnarray*}
\end{thm}
Here it is understood that
\begin{eqnarray}  \label{understood}
 \frac1{1-uq^m} =\left\{
\begin{array}{ll}
\sum_{k=0}^\infty (uq^m)^k, & \text{ if }\; m \ge 0, \\
-\sum_{k=1}^\infty (u^{-1}q^{-m})^k, & \text{ if }\; m < 0.
\end{array}
\right.
\end{eqnarray}

\begin{proof} Using Proposition \ref{111} (in the fourth line
below) and (\ref{exponential}), we have
{\allowdisplaybreaks\begin{eqnarray*}
 && \frac1{(u)_{\infty}(u^{-1}q)_{\infty}} \\
  &=& \sum_{l=0}^\infty \frac{u^l}{(q)_l} \sum_{n=0}^\infty
  \frac{(u^{-1}q)^n}{(q)_n} \\
  &=& \sum_{k\ge 0}  u^k \sum_{n=0}^\infty \frac{q^n}{(q)_n (q)_{n+k}}
   + \sum_{k >0} u^{-k} \sum_{n=0}^\infty \frac{q^{n+k}}{(q)_n (q)_{n+k}} \\
  &=& \frac1{(q)_{\infty}^2}
  (\sum_{k\ge 0}   \sum_{m=0}^\infty (-1)^m q^{m(m+1)/2} q^{km} u^k
   + \sum_{k >0}  \sum_{m=0}^\infty (-1)^m q^{m(m+1)/2} q^{km}q^k u^{-k}) \\
  &=& \frac1{(q)_{\infty}^2}   \sum_{m=0}^\infty (-1)^m q^{\hf m(m+1)}
   \left( \sum_{k\ge 0}   q^{km} u^k
   + \sum_{k >0} q^{k(m+1)} u^{-k} \right) \\
   &=&  \frac1{(q)_{\infty}^2} \sum_{m \in \Z} (-1)^m q^{\hf m(m+1)} \frac1{1-uq^m}.
\end{eqnarray*}}
The last equation follows by the interpretation
(\ref{understood}).
\end{proof}

\subsection{Lie algebra $\ainf$}

Denote by $\gl$ the Lie algebra of all matrices $(a_{ij})_{i,j \in
\Z}$ satisfying $a_{ij} =0$ for $|i -j|\gg 0$. Denote by $E_{ij}$
the infinite matrix with $1$ at $(i, j)^{th}$ place and $0$
elsewhere and let the weight of $ E_{ij}$ be $j - i $. This defines
a $\Z$--principal gradation $\gl = \bigoplus_{j \in \Z} \gl_j$.
Denote by $\ainf \equiv \ohgl_\infty = \gl \oplus \C C$ the central
extension given by the following $2$--cocycle with values in $\C$
(cf. \cite{DJKM1}):
\begin{eqnarray}
 \alpha(A, B) = \trace \left([J, A]B \right),
  \label{eq_cocy}
\end{eqnarray}
where $J = \sum_{ j \leq 0} E_{jj}$. The $\Z$--gradation of the Lie
algebra $\gl$ extends to $\hgl$ by letting the weight of $C$ be $0$.
This leads to a triangular decomposition (i.e. a direct sum of
subspaces of positive, zero, and negative weights):
$$\hgl = (\ainf)_{+} \oplus (\ainf)_{0} \oplus (\ainf)_{-}, $$
where $(\ainf)_{0} = \gl_0 \oplus \C C.$
Let
\begin{eqnarray*}
 H^a_i  = E_{ii} - E_{i+1, i+1} + \delta_{i,0} C, \quad i \in \Z.
\end{eqnarray*}
Denote by $L(\hgl; \Lambda,k)$ the highest weight $\hgl$--module
with highest weight $\Lambda  \in (\ainf)_{0}^{*}$ and level $k$,
where $C$ acts as a scalar $k\cdot I$. Let $\hL_j^a \in
(\ainf)_0^{*}$ be the fundamental weights, i.e. $\hL_j^a ( H_i^a ) =
\delta_{ij}.$ The Dynkin diagram for $\hgl$, with fundamental
weights labeled, is the following:
  \begin{equation*}
  \begin{picture}(150,45) 
  \put(-16,23){\dots}
  \put(0,20){$\circ$}
  \put(7,23){\line(1,0){32}}
  \put(40,20){$\circ$}
  \put(47,23){\line(1,0){32}}
  \put(81,20){$\circ$}
  \put(88,23){\line(1,0){32}}
  \put(121,20){$\circ$}
  \put(128,23){\line(1,0){32}}
  \put(161,20){$\circ$}
  \put(170,23){\dots}
  \put(-4,9){$-2$}
  \put(36,9){$-1$}
  \put(81,9){$0$}
  \put(121,9){$1$}
  \put(162,9){$2$}
  \end{picture}
  \end{equation*}

\subsection{The $1$-point $\hgl$-functions of
level $-1$}

Consider a pair of free bosonic ghosts
\begin{equation*}
\gamma^{\pm}(z)=\sum_{r\in\hf+\Z}\gamma^\pm_r z^{-r-\hf},
\end{equation*}
with nontrivial commutation relations
$$
 [\gamma^+_r,\gamma^-_s]=\delta_{r+s,0},
 \quad r, s \in \hf+\Z.
$$
Let $\F^{-1}$ denote the Fock space generated by $\gamma^\pm(z)$
with vacuum vector $\vac$ (cf. \cite{K1} for more on Fock spaces
and normal ordering $\no \; \no$ in vertex algebras). An action of
$\hgl$ of level $-1$ on $\F^{-1}$ is given by (cf. e.g.
\cite{W1}):
\begin{equation*}
E(z,w) \equiv
\sum_{i,j\in\Z}E_{ij}z^{i-1}w^{-j}=-\no\gamma^+(z)\gamma^-(w)\no,
\end{equation*}
so that $E_{ij}=-\no\gamma^+_{-i+\hf}\gamma^-_{j-\hf}\no$. The
Virasoro field is given by
\begin{equation*}
L(z)=\sum_{n\in\Z}L_nz^{-n-2} =\hf\big{(}\no\gamma^+(z)\partial
\gamma^-(z)\no -\no\partial\gamma^+(z)\gamma^-(z)\no\big{)},
\end{equation*}
and we have $[L_0,\gamma^{\pm}_r]=-r\gamma^\pm_r$. According to
the eigenvalues of the charge operator
$e_{11}=\sum_{r\in\hf+\Z}\no\gamma^+_r\gamma^-_{-r}\no$, the
$\hgl$-module $\F^{-1}$ has the following decomposition:
$$\F^{-1}=\bigoplus_{n\in\Z}\F^{-1}_{(n)}.$$

Following Bloch-Okounkov \cite{BO}, we introduce the following
operators in $\ainf$:
$$\begin{aligned}
\no\A(t)\no
&=\sum_{k\in\Z}E_{kk}t^{k-\hf},\\
\A(t)&=\no\A(t)\no+\frac{C}{t^\hf-t^{-\hf}}.
\end{aligned}$$
When acting on $\F^{-1}$, $\A(t) =
-\sum_{r\in\hf+\Z}t^r\gamma^+_{-r}\gamma^-_{r}$.

Define the Bloch-Okounkov $n$-point $\ainf$-correlation function
(or $n$-point $\ainf$-function for short) of level $-1$
(associated to $m\in \Z$) to be
\begin{equation*}
\mathfrak{A}^{(m)}_{-1}(q; t_1,\ldots,t_n)
:=\trace_{\F^{-1}_{(m)}}(q^{L_0}\A(t_1) \cdots \A(t_n)).
\end{equation*}

For a partition $\la$ we use $\ell(\la)$ and $|\la|$ to denote the
length and the size of $\la$ respectively.

\begin{lem}\label{222} We have
\begin{itemize}
\item[(i)]
$\displaystyle\sum_{\ell(\la)=l}q^{|\la|}=\frac{q^l}{(q)_l}$,
\item[(ii)] $\displaystyle\sum_{\ell(\la)=l}q^{|\la|}t^{\la_i}
=\frac{tq^l}{(1-q)\cdots(1-q^{i-1})(1-q^it)\cdots(1-q^lt)}$,\quad
$i\le l$.
\end{itemize}
\end{lem}

\begin{proof}
Part (i) follows from the identity $\sum_{\ell(\la)\le
l}q^{|\la|}= (q)_l^{-1}$, while (ii) follows from another
well-known identity
$$\sum_{\ell(\la)\le
l}q^{|\la|}t^{\la_i}=\displaystyle\frac{1}{(1-q)\cdots(1-q^{i-1})(1-q^it)\cdots(1-q^lt)}.$$
\end{proof}

\begin{thm}
The one-point function $\mathfrak{A}^{(0)}_{-1}(q; t)$ is given by
\begin{align*}
\frac{_2\Phi_1(0,0;q;q)}{t^{-\hf}-t^{\hf}}+t^{\hf}\sum_{i=1}^\infty\frac{q^{i-1}}{(q)_{i-1}(q)_{i-1}}\big{(}
{}_3\Phi_2(0,0,q;tq^i,q^i;q)-1\big{)}\\
-t^{-\hf}\sum_{i=1}^\infty\frac{q^{i-1}}{(q)_{i-1}(q)_{i-1}}\big{(}
{}_3\Phi_2(0,0,q;t^{-1}q^i,q^i;q)-1\big{)}.
\end{align*}
\end{thm}

\begin{proof}
Note that $\gamma^+_{-r_1}\gamma^+_{-r_2}\cdots \gamma^+_{-r_l}
\gamma^-_{-s_1}\cdots\gamma^-_{-s_l}\vac$ is an eigenvector of
$q^{L_0}\A(t)$ of eigenvalue
$q^{\sum_{i=1}^l (r_i+s_i)}
\left(\sum_{i=1}^l(t^{r_i}-t^{-s_i})-\sum_{r\in\hf+\Z_+}t^{-r}
\right). $ Since $\F^{-1}_{(0)}$ has a basis given by
$\gamma^+_{-r_1} \cdots\gamma^+_{-r_l} \gamma^-_{-s_1}
\cdots\gamma^-_{-s_l}\vac$ for $r_1\ge \cdots\ge
r_l>0;s_1\ge\cdots\ge s_l>0, l \ge 0$, we have by Lemma~\ref{222}
that
\begin{align}
\mathfrak{A}^{(0)}_{-1}& (q;t) \label{onepoint} \\
 & =
\sum_{l=0}^\infty q^{-l} \sum_{\ell(\la)=\ell(\mu)=l}
q^{|\la|+|\mu|}\left(\sum_{i=1}^l(t^{\la_i-\hf}-t^{-\mu_i+\hf}) +
\frac{1}{t^{-\hf}-t^{\hf}} \right) \nonumber \\
 & =
\displaystyle\frac{_2\Phi_1(0,0;q;q)}{t^{-\hf}-t^{\hf}}+\sum_{l=1}^\infty
q^{-l} \sum_{\ell(\la)=\ell(\mu)=l}
q^{|\la|+|\mu|}\big{(}\sum_{i=1}^l(t^{\la_i-\hf}-t^{-\mu_i+\hf})
\big{)}.  \nonumber
\end{align}

We compute that
{\allowdisplaybreaks
\begin{align*}
&\sum_{l=1}^\infty q^{-l} \sum_{\ell(\la)=\ell(\mu)=l}
q^{|\la|+|\mu|}\sum_{i=1}^l t^{\la_i-\hf} \\
& =t^{-\hf}\sum_{l=1}^\infty q^{-l}\sum_{i=1}^l
\frac{q^{2l}t}{(1-q)\cdots(1-q^{i-1})(1-q^it)\cdots(1-q^lt)(q)_l} \\
& = t^{\hf}\sum_{i=1}^\infty\sum_{l=i}^\infty
\frac{q^l}{(q)_{i-1}(q)_{l}(q^it)_{l-i+1}} \\
& =t^{\hf}\sum_{i=1}^\infty\frac{q^{i-1}}{(q)_{i-1}(q)_{i-1}}
\sum_{s=1}^\infty\frac{q^s}{(tq^i)_s(q^i)_s} \\
&=t^{\hf}\sum_{i=1}^\infty\frac{q^{i-1}}{(q)_{i-1}(q)_{i-1}}
\big{(} {}_3\Phi_2(0,0,q;tq^i,q^i;q)-1\big{)}.
\end{align*}}
Now the theorem follows from this computation and
(\ref{onepoint}).
\end{proof}
\subsection{The generalized 1-point $\hgl$-function} \label{gl:subsection:1to2}

Let $A$ be the operator on $\F^{-1}$ acting trivially on $\vac$
such that $[A,\gamma^{+}_r]=\gamma^+_r$, and $[A,\gamma^-_r]=0$,
for all $r\in\hf+\Z$.  Let $B$ be the operator acting trivially on
$\vac$ with $[B,\gamma^{-}_r]=\gamma^-_r$, and $[B,\gamma^+_r]=0$,
for all $r\in\hf+\Z$.

\begin{thm} \label{gl:general:AB}
We have
{\allowdisplaybreaks\begin{align*}
\trace_{\F^{-1}}&(q^{L_0}x^Ay^B\A(t))=
\frac{1}{(t^{-\hf}-t^\hf)(xq^\hf)_\infty(yq^\hf)_\infty}\\
&+\frac{x(tq)^\hf(xtq^{\frac{3}{2}})_\infty}{(1-xq^\hf)(tq)_\infty(xq^\hf)_\infty(yq^\hf)_\infty}
{}_2\Phi_2(xq^\hf,xq^\hf;xq^{\frac{3}{2}},txq^{\frac{3}{2}};tq^2)\\
&-\frac{y(t^{-1}q)^\hf(yt^{-1}q^{\frac{3}{2}})_\infty}{(1-yq^\hf)(t^{-1}q)_\infty(xq^\hf)_\infty(yq^\hf)_\infty}
{}_2\Phi_2(yq^\hf,yq^\hf;yq^{\frac{3}{2}},yt^{-1}q^{\frac{3}{2}};t^{-1}q^2).
\end{align*}}
\end{thm}

\begin{proof}

By Lemma \ref{222} and (\ref{exponential}), we have
\begin{align*}
\sum_{k=0}^\infty
q^{-\frac{k}{2}}z^k\sum_{\ell(\mu)=k}q^{|\mu|}=\frac{1}{(zq^\hf)_\infty}.
\end{align*}
Since the vector
$\gamma^+_{-r_1}\gamma^+_{-r_2}\cdots\gamma^+_{-r_l}
\gamma^-_{-s_1}\cdots\gamma^-_{-s_k}\vac$
is an eigenvector of $x^Ay^Bq^{L_0} \A(t)$ of eigenvalue
$x^ly^kq^{\sum_{i=1}^l r_i+\sum_{j=1}^ks_j}
\big{(}\sum_{i=1}^lt^{r_i}-\sum_{j=1}^kt^{-s_j}-\sum_{r\in\hf+\Z_+}t^{-r}
\big{)},$ we have
\begin{align*}
\trace_{\F^{-1}} &(x^Ay^Bq^{L_0}\A(t)) \\
 =& \sum_{k=0}^\infty
y^kq^{-\frac{k}{2}} \sum_{\ell(\mu)=k}q^{|\mu|} \sum_{l=0}^\infty
x^lq^{-\frac{l}{2}} \sum_{\ell(\la)=l}q^{|\la|}
\big{(}\sum_{i=1}^lt^{\la_i-\hf}-\sum_{j=1}^kt^{-\mu_j+\hf}\big{)}
\\
& +\frac{1}{(xq^\hf)_\infty(yq^\hf)_\infty(t^{-\hf}-t^{\hf})}.
\end{align*}
Next we claim that
\begin{align}\label{555}
\sum_{l=1}^\infty & x^lq^{-\frac{l}{2}}\sum_{\ell(\la)
=l}q^{|\la|}\sum_{i=1}^lt^{\la_i-\hf} \\
 & =
\frac{x(tq)^\hf(xtq^{\frac{3}{2}})_\infty}{(1-xq^\hf)(tq)_\infty(xq^\hf)_\infty}
{}_2\Phi_2(xq^\hf,xq^\hf;xq^{\frac{3}{2}},xtq^{\frac{3}{2}};tq^2).
\nonumber
\end{align}
To see this we compute directly
{\allowdisplaybreaks
\begin{align*}
&\sum_{l=1}^\infty
x^lq^{-\frac{l}{2}}\sum_{\ell(\la)=l}q^{|\la|}\sum_{i=1}^lt^{\la_i-\hf}\\
&= t^\hf \sum_{l=1}^\infty \sum_{i=1}^\infty
\frac{x^lq^{\frac{l}{2}}}{(q)_{i-1}(1-tq^i)\cdots(1-tq^l)}\\
&=t^\hf\sum_{i=1}^\infty\frac{1}{(q)_{i-1}(tq^i)_\infty}\sum_{l=i}^\infty
x^lq^{\frac{l}{2}}(tq^{l+1})_\infty\\
&=t^\hf\sum_{i=1}^\infty\frac{1}{(q)_{i-1}(tq^i)_\infty}\sum_{l=i}^\infty
x^lq^{\frac{l}{2}}\sum_{m=0}^\infty \frac{q^{\hf
m(m-1)}q^{m(l+1)}(-t)^m}{(q)_m}\\
&=t^\hf\sum_{m=0}^\infty \frac{q^{\hf
m(m-1)}(-tq)^m}{(q)_m}\sum_{i=1}^\infty\frac{1}{(q)_{i-1}(tq^i)_\infty}\sum_{l=i}^\infty
x^lq^{l(m+\hf)}\\
&=\frac{t^\hf}{(tq)_\infty}\sum_{m=0}^\infty\frac{q^{\hf
m(m-1)}(-tq)^m}{(q)_m}\sum_{i=1}^\infty\frac{(tq)_{i-1}}{(q)_{i-1}}
\frac{x^iq^{i(m+\hf)}}{1-xq^{m+\hf}}\\
&=\frac{xt^\hf}{(tq)_\infty} \sum_{m=0}^\infty\frac{q^{\hf
m(m-1)}(-tq)^m q^{m+\hf}}{(q)_m(1-xq^{m+\hf})}\sum_{i=1}^\infty
\frac{(tq)_{i-1}}{(q)_{i-1}} (xq^{(m+\hf)})^{i-1}\\
&=\frac{xt^\hf}{(1-xq^\hf)(tq)_\infty} \sum_{m=0}^\infty
\frac{q^{\hf m(m-1)}(-tq)^mq^{m+\hf}(xq^\hf)_m}{(q)_m
(xq^{\frac{3}{2}})_m}
\frac{(xtq^{m+\frac{3}{2}})_\infty}{(xq^{m+\hf})_\infty}\\
&=\frac{x(tq)^\hf(xtq^{\frac{3}{2}})_\infty}{(1-xq^\hf)(tq)_\infty(xq^\hf)_\infty}
\sum_{m=0}^\infty \frac{q^{\hf
m(m-1)}(-tq^2)^m(xq^\hf)^2_m}{(q)_m(xq^\frac{3}{2})_m(xtq^{\frac{3}{2}})_m}\\
&=
\frac{x(tq)^\hf(xtq^{\frac{3}{2}})_\infty}{(1-xq^\hf)(tq)_\infty(xq^\hf)_\infty}
{}_2\Phi_2(xq^\hf,xq^\hf;xq^{\frac{3}{2}},xtq^{\frac{3}{2}};tq^2).
\end{align*} }

The theorem now follows from this and a similar expression for the
other summation (involving $y^B$).
\end{proof}

\begin{rem} \label{rem:residue}
Noting that the operator $B-A$ is the same as the operator
$e_{11}$ defined earlier, we have
$$\mathfrak{A}^{(m)}_{-1}(q; t) =[z^m]
\trace_{\F^{-1}}(q^{L_0}x^Ay^B\A(t))|_{x=z^{-1},y=z}.
$$
Here and below $[z^m]X$ denotes the coefficient of $z^m$ in the
expansion of $X$ in $z$. For this reason, we refer to
$\trace_{\F^{-1}} (q^{L_0}x^Ay^B\A(t))$ as the generalized 1-point
$\ainf$-function.
\end{rem}

\subsection{The generalized $2$-point $\ainf$-function}

We will compute the generalized $2$-point $\ainf$-function
$\trace_{\F^{-1}}q^{L_0} x^A y^B \A (t_1)\A (t_2)$. Similar to the
computation of the generalized $1$-point function in the previous
subsection, the calculation of the generalized 2-point function
essentially boils down to the following calculation:

{\allowdisplaybreaks
\begin{align*}
&\sum_{i=1}^\infty\sum_{j=i+1}^\infty\sum_{l=j}^\infty \frac
{x^lq^{\frac{l}{2}}t_1t_2}
{(q)_{i-1}(1-q^it_1)\cdots(1-q^{j-1}t_1)(1-q^jt_1t_2)\cdots
(1-q^lt_1t_2)}\\
&=\frac{t_1t_2}{(qt_1t_2)_\infty}\sum_{i=1}^\infty\frac{1}{(q)_{i-1}}
\sum_{j=i+1}^\infty
\frac{(qt_1t_2)_{j-1}}{(1-q^it_1)\cdots(1-q^{j-1}t_1)}\sum_{l=j}^\infty
x^lq^{\frac{l}{2}}(q^{l+1}t_1t_2)_\infty \\
&=\frac{t_1t_2}{(qt_1)_\infty(qt_1t_2)_\infty}\sum_{i=1}^\infty\frac{(qt_1)_{i-1}}{(q)_{i-1}}
\sum_{j=i+1}^\infty (q^jt_1)_\infty (qt_1t_2)_{j-1}\sum_{l=j}^\infty
x^lq^{\frac{l}{2}}(q^{l+1}t_1t_2)_\infty\\
&=\frac{t_1t_2}{(qt_1)_\infty(qt_1t_2)_\infty}\sum_{i=1}^\infty\frac{(qt_1)_{i-1}}{(q)_{i-1}}
\sum_{j=i+1}^\infty (q^jt_1)_\infty (qt_1t_2)_{j-1}
\\
&\qquad\qquad\qquad\qquad\times\sum_{s=0}^\infty
\frac{(-qt_1t_2)^s q^{\hf s(s-1)}}{(q)_s}\sum_{l=j}^\infty
(xq^{s+\hf})^l\\
&=\frac{t_1t_2}{(qt_1)_\infty(qt_1t_2)_\infty}
\\
& \quad \times \sum_{i=1}^\infty\frac{(qt_1)_{i-1}}{(q)_{i-1}}
\sum_{j=i+1}^\infty (q^jt_1)_\infty
(qt_1t_2)_{j-1}\sum_{s=0}^\infty
\frac{(-qt_1t_2)^s q^{\hf s(s-1)}(xq^{s+\hf})^j}{(q)_s(1-xq^{s+\hf})}\\
&=\frac{t_1t_2}{(qt_1)_\infty(qt_1t_2)_\infty}
\\
& \quad \times \sum_{i=1}^\infty \frac{(qt_1)_{i-1}}{(q)_{i-1}}
\sum_{s=0}^\infty \frac{(-qt_1t_2)^s q^{\hf
s(s-1)}}{(q)_s(1-xq^{s+\hf})} \sum_{j=i+1}^\infty
(xq^{s+\hf})^j(q^jt_1)_\infty (qt_1t_2)_{j-1},\end{align*}} which
is equal to {\allowdisplaybreaks
\begin{align*}
&\frac{t_1t_2}{(qt_1)_\infty}\sum_{i=1}^\infty\frac{(qt_1)_{i-1}}{(q)_{i-1}}
\sum_{s=0}^\infty \frac{(-qt_1t_2)^s q^{\hf
s(s-1)}}{(q)_s(1-xq^{s+\hf})}
\sum_{j=i+1}^\infty\frac{(q^jt_1)_\infty}{(q^jt_1t_2)_\infty} (xq^{s+\hf})^j\\
&=\frac{t_1t_2}{(qt_1)_\infty}\sum_{i=1}^\infty\frac{(qt_1)_{i-1}}{(q)_{i-1}}
\sum_{s=0}^\infty \frac{(-qt_1t_2)^s q^{\hf
s(s-1)}}{(q)_s(1-xq^{s+\hf})}
\sum_{j=i+1}^\infty\sum_{m=0}^\infty\frac{(t_2^{-1})_m}{(q)_m} (q^jt_1t_2)^m(xq^{s+\hf})^j\\
&=\frac{t_1t_2}{(qt_1)_\infty}\sum_{i=1}^\infty\frac{(qt_1)_{i-1}}{(q)_{i-1}}
\sum_{s=0}^\infty \frac{(-qt_1t_2)^s q^{\hf
s(s-1)}}{(q)_s(1-xq^{s+\hf})}
\sum_{m=0}^\infty\frac{(t_2^{-1})_m}{(q)_m}(t_1t_2)^m\sum_{j=i+1}^\infty (xq^{m+s+\hf})^j\\
&=\frac{t_1t_2}{(qt_1)_\infty}\sum_{i=1}^\infty\frac{(qt_1)_{i-1}}{(q)_{i-1}}
\sum_{s=0}^\infty \frac{(-qt_1t_2)^s q^{\hf
s(s-1)}}{(q)_s(1-xq^{s+\hf})}
\sum_{m=0}^\infty\frac{(t_2^{-1})_m}{(q)_m}(t_1t_2)^m\frac{(xq^{m+s+\hf})^{i+1}}{1-xq^{m+s+\hf}}\\
&=\frac{t_1t_2}{(qt_1)_\infty} \sum_{s=0}^\infty
\frac{(-qt_1t_2)^s
q^{\hf s(s-1)}}{(q)_s(1-xq^{s+\hf})}\\
&\qquad\qquad\times\sum_{m=0}^\infty\frac{(t_2^{-1})_m}{(q)_m}(t_1t_2)^m\frac{(xq^{m+s+\hf})^{2}}{1-xq^{m+s+\hf}}
\sum_{i=1}^\infty\frac{(qt_1)_{i-1}}{(q)_{i-1}} (xq^{m+s+\hf})^{i-1}\\
&=\frac{t_1t_2}{(qt_1)_\infty} \sum_{s=0}^\infty \frac{(-qt_1t_2)^s
q^{\hf s(s-1)}}{(q)_s(1-xq^{s+\hf})}
\sum_{m=0}^\infty\frac{(t_2^{-1})_m}{(q)_m}(t_1t_2)^m\frac{(xq^{m+s+\hf})^{2}}{1-xq^{m+s+\hf}}
\frac{(xq^{m+s+\frac{3}{2}}t_1)_\infty}{(xq^{m+s+\hf})_\infty}\\
&=\frac{x^2qt_1t_2(xq^{\frac{3}{2}}t_1)_\infty}{(1-xq^\hf)^2(xq^\hf)_\infty(qt_1)_\infty}\\
&\qquad\qquad\times \sum_{s=0}^\infty
\frac{(xq^\hf)^3_s(-q^3t_1t_2)^s q^{\hf
s(s-1)}}{(xq^{\frac{3}{2}}t_1)_s(q)_s(xq^{\frac{3}{2}})_s^2}
\sum_{m=0}^\infty\frac{(t_2^{-1})_m(q^2t_1t_2)^m(xq^{s+\hf})^2_m}
{(q)_m(xq^{s+\frac{3}{2}})_m(xq^{s+\frac{3}{2}}t_1)_m}\\
&=\frac{x^2qt_1t_2(xq^{\frac{3}{2}}t_1)_\infty}{(1-xq^\hf)^2(xq^\hf)_\infty(qt_1)_\infty}
\\
&\times\sum_{s=0}^\infty \frac{(xq^\hf)^3_s(-q^3t_1t_2)^s q^{\hf
s(s-1)}}{(xq^{\frac{3}{2}}t_1)_s(q)_s(xq^{\frac{3}{2}})_s^2}
{}_3\Phi_2
(t_2^{-1},xq^{s+\hf},xq^{s+\hf};xq^{s+\frac{3}{2}},
xq^{s+\frac{3}{2}}t_1;q^2t_1t_2).
\end{align*}}
Denote the above expression by $\overline{\Gamma}(x,t_1,t_2)$. We
point out that all the 3-2 $q$-hypergeometric series in
$\overline{\Gamma}(x,t_1,t_2)$ are of type II, meaning that
${}_3\Phi_2(a,b,c;d,e;z)$ satisfies $\frac{de}{abc}=z$.  Define
\begin{equation*}
\Gamma(x,y,t_1,t_2):=\frac{(t_1t_2)^{-\hf}}{(q^\hf
y)_\infty}\left(
\overline{\Gamma}(x,t_1,t_2)+\overline{\Gamma}(x,t_2,t_1)\right).
\end{equation*}
To simplify notation we set
\begin{equation*}
\Omega(x,y,t):=
\frac{x(tq)^\hf(xtq^{\frac{3}{2}})_\infty}{(1-xq^\hf)(tq)_\infty(xq^\hf)_\infty(yq^\hf)_\infty}
{}_2\Phi_2(xq^\hf,xq^\hf;xq^{\frac{3}{2}},txq^{\frac{3}{2}};tq^2).
\end{equation*}

The following is now a straightforward calculation, based on the
above calculations.

\begin{thm}
The generalized $2$-point $\ainf$-function
$\trace_{\F^{-1}}q^{L_0} x^A y^B \A (t_1)\A (t_2)$ is equal to
{\allowdisplaybreaks
\begin{align*}
&\Gamma(x,y,t_1,t_2)+\Gamma(y,x,t^{-1}_1,t^{-1}_2)
+\Omega(x,y,t_1t_2)+\Omega(y,x,t^{-1}_1t_2^{-1})\\
&+\frac{1}{t_1^{-\hf}-t_1^\hf}\left(
\Omega(x,y,t_2)-\Omega(y,x,t_2^{-1}) \right)
+\frac{1}{t_2^{-\hf}-t_2^\hf}\left(
\Omega(x,y,t_1)-\Omega(y,x,t_1^{-1})
\right)\\
&-(xq^\hf)_\infty(yq^\hf)_\infty\left(\Omega(x,y,t_1)\Omega(y,x,t_2^{-1})
+\Omega(x,y,t^{-1}_1)\Omega(y,x,t_2)\right)\\
&+\frac{1}{(t_1^{-\hf}-t_1^\hf)(t_2^{-\hf}-t_2^\hf)(xq^\hf)_\infty(yq^\hf)_\infty}.
\end{align*}}
\end{thm}

\begin{rem}  \label{rem:residue2}
The $2$-point $\ainf$-function can be recovered from the
generalized  $2$-point $\ainf$-function:
$$
\mathfrak{A}^{(m)}_{-1}(q;t_1,t_2) =[z^m]
\trace_{\F^{-1}}(q^{L_0}x^Ay^B\A(t_1) \A(t_2))|_{x=z^{-1},y=z}.
$$
\end{rem}

\subsection{The $n$-point $\hgl$-functions of level $-l$}

Let $l\in \N$. Recall a generalized partition $\la$ of depth $l$
is an ordered $l$-tuple of non-increasing integers; that is, $\la
= (\la_1,\la_2,\dots,\la_l)$ where $\la_1\geq \cdots\geq\la_l$.
Associated to a generalized partition $\la =
(\la_1,\la_2,\dots,\la_l)$ we define a highest weight
$\Lambda(\la)$ of $\hgl$ via
$$\Lambda(\la) = (\la_l-\la_1-l)\Lambda_0^a +
\sum_{k=1}^{i-1}(\la_k-\la_{k+1})\Lambda_l^a + \la_i\Lambda_i^a,$$
where $i$ is the index of the last positive entry of $\la$.

We let $\F^{-l}:=(\F^{-1})^{\otimes l}$ to denote the Fock space
generated by $l$ pairs of free bosonic fields
$\gamma^{\pm,i}(z)=\sum_{r\in\frac{1}{2}+\Z}\gamma^{\pm,i}_rz^{-r-1/2}$,
$i=1,\cdots,l$, with non-trivial commutation relations
$[\gamma^{+,i}_r,\gamma^{-,j}_s]=\delta_{ij}\delta_{r+s,0}$. On
$\F^{-l}$ there is an action of $a_\infty$ of level $-l$ given by
\begin{equation*}
E(z,w)=-\sum_{i=1}^l\no\gamma^{i,+}(z)\gamma^{i,-}(w)\no.
\end{equation*}
Furthermore there exists a commuting action of the general linear
algebra $\gl(l)$ whose elementary matrices $e_{ij}$ acts by the
formula
\begin{equation*}
e_{ij}=- \sum_{r\in1/2+\Z}:\gamma^{+,i}_{-r}\gamma^{-,j}_{r}:.
\end{equation*}
This action lifts to that of $GL(l)$.

\begin{prop} \label{gl:duality} \cite{KR} (also see \cite[Theorem
5.1]{W1})
We have the following $(GL(l),\hgl)$-module decomposition:
$$\F^{-l} \cong \bigoplus_{\la} V_\la(GL(l))\otimes
L(\hgl;\Lambda(\la),-l),
$$
where $\la$ is a generalized partition of depth $l$ and
$V_\la(GL(l))$ is the irreducible $GL(l)$-module of highest weight
$\la$.
\end{prop}
The Bloch-Okounkov $n$-point $\ainf$-function of level $-l$
(associated to a generalized partition $\la$ of depth $l$) is
defined as
$$\mathfrak A^\la_{-l}(q; t_1,\dots,t_n) =
\trace_{L(\hgl;\Lambda(\la),-l)}q^{L_0}\A(t_1)\A(t_2)\cdots\A(t_n).
$$

\begin{thm} \label{gl:npt:lvll}
The $n$-point $\hgl$-function of level $-l$,
$\mathfrak{A}^\la_{-l}(q;t_1,\dots,t_n)$, is equal to
$$
\sum_{\sigma\in S_l} (-1)^{\ell(\sigma)}\mathfrak
A^{(k_1)}_{-1}(q;t_1,\dots,t_n)\cdots \mathfrak
A^{(k_l)}_{-1}(q;t_1,\dots,t_n),
$$
where $k_i = (\la+\rho-\sigma(\rho),\varepsilon_i)$.
\end{thm}

\begin{proof}
We shall denote by $e_{ij}, 1\le i,j \le l,$ the standard matrix
elements in the Lie algebra $\mathfrak{gl}(l)$. Applying
$\trace_{\F^{-l}} z_1^{e_{11}}\cdots z_l^{e_{ll}} q^{L_0}
\A(t_1)\cdots \A(t_n)$ to both sides of the identity in
Proposition~\ref{gl:duality}, we obtain that
$$
\prod_{i=1}^l \trace_{\F^{-1}} z_i^{e_{ii}} q^{L_0} \A(t_1)\cdots
\A(t_n) = \sum_\la \mathsf{ch}_\la^{gl}(z_1,\dots,z_l)\mathfrak
A^\la_{-l}(q;t_1,\dots,t_n),
$$
since $z_i^{e_{ii}}$ acts on the left-hand side only on the
$i^{th}$ tensor factor. On the other hand, we have
$$\begin{aligned}
\prod_{i=1}^l\trace_{\F^{-1}} z_i^{e_{ii}} q^{L_0} &\A(t_1)\cdots
\A(t_n)
= \prod_{i=1}^l \left(\sum_{m_i\in\Z} z_i^{m_i}\mathfrak
A^{(m_i)}_{-1}(q;t_1\dots,t_n)\right) \\ &= \sum_{\mathbf{m}\in\Z^l}
z_1^{m_1}\cdots z_l^{m_l}\mathfrak
A^{(m_1)}_{-1}(q;t_1,\dots,t_n)\cdots\mathfrak
A^{(m_l)}_{-1}(q;t_1,\dots,t_n).\end{aligned}$$

Recall (cf. \cite[pp.~399]{FH}) that
\begin{eqnarray}  \label{weylchar}
 \mathsf{ch}_\la^{gl}(z_1,\dots,z_l) =
\frac{|z_j^{\la_i+l-i}|}{|z_j^{l-i}|},
\end{eqnarray}
where $|a_{ij}|$ denotes the determinant of the matrix $(a_{ij})$.

Combining the above identities with the Weyl denominator formula,
we obtain that
\begin{align*}
\sum_{\la}  &
 |z_j^{\la_i+l-i}|\, \mathfrak A^\la_{-l}(q;t_1,\dots,t_n) \\
 &= \sum_{\sigma\in S_l}
(-1)^{\ell(\sigma)} \mathbf{z}^{\sigma(\rho)}\prod_{i=1}^l
\trace_{\F^{-1}} z_i^{e_{ii}} q^{L_0} \A(t_1)\cdots \A(t_n) \\
 &= \sum_{\sigma\in S_l}
(-1)^{\ell(\sigma)} \mathbf{z}^{\sigma(\rho)}
\sum_{\mathbf{m}\in\Z^l} z_1^{m_1}\cdots z_l^{m_l}\mathfrak
A^{(m_1)}_{-1}(q;t_1,\dots,t_n)\cdots\mathfrak
A^{(m_l)}_{-1}(q;t_1,\dots,t_n),
\end{align*}
where we denote $\mathbf{z}^{\sigma(\rho)}=\prod_{i=1}^l
z_i^{(\sigma(\rho),\varepsilon_i)}$.

Now the theorem follows from comparing the coefficients of
$\prod_{i=1}^l z_i^{\la_i+l-i}$ on both sides of the above
identity.
\end{proof}

\begin{rem}
Combining Theorem~\ref{gl:npt:lvll} with the calculation of the
$1$-point and $2$-point functions of level $-1$ in the previous
subsections, we have computed the $1$-point and $2$-point
$\ainf$-functions of level $-l$.
\end{rem}

\subsection{A $q$-difference equation for $\ainf$-functions}

Even though it is difficult to compute the $n$-point function in
general, we have the following $q$-difference equation of level
$-1$ which could be helpful (compare \cite{BO, Ok} for a
difference equation of level $1$).

\begin{thm}
The $n$-point $\ainf$-function of level $-1$ satisfies the
following $q$-difference equation:
\begin{align*}
\mathfrak A^{(0)}_{-1} & (q; qt_1,\dots,t_n) \\
& = \sum_{s=0}^{n-1} (-1)^{s+1} \sum_{1<i_1<\cdots<i_s\leq n}
\mathfrak A^{(0)}_{-1}(q; t_1t_{i_1}\cdots
t_{i_s},\dots,\widehat{t}_{i_1},\dots,\widehat{t}_{i_s},\dots),
\end{align*}
where $\,\widehat{}\,$ denotes a deleted term.
\end{thm}

\begin{proof}
Applying the commutation relation for $\A(t)$ and $\gamma_r^-$
repeatedly yields $$\A(t_2)\cdots\A(t_n)\gamma_r^- =
\sum_{P\subset \{2,\dots,n\}} (-1)^{|P|} \left(\prod_{i\in P}
t_i^r\right)\gamma_r^-\prod_{i\not\in P} \A(t_i).
$$
Below we write $\trace =\trace_{\F_{(0)}^{-1}}$. Taking trace of
the above gives
$$
\trace q^{L_0}\gamma_{-r}^+\A(t_2)\cdots\A(t_n)\gamma_r^- =
\sum_{P\subset \{2,\dots,n\}} (-1)^{|P|} \left(\prod_{i\in P}
t_i^r\right)\trace q^{L_0}\gamma_{-r}^+\gamma_r^-\prod_{i\not\in
P} \A(t_i).
$$
It follows from the commutation relation for $L_0$ that
$$
\begin{aligned}
-\trace q^{L_0}\gamma_{-r}^+\A(t_2)\cdots\A(t_n)\gamma_r^- &=
-\trace \gamma_r^-q^{L_0}\gamma_{-r}^+\A(t_2)\cdots\A(t_n) \\ &=
-q^r\trace
q^{L_0}\gamma_r^-\gamma_{-r}^+\A(t_2)\cdots\A(t_n).
\end{aligned}
$$

Multiplying both sides of the above by $t_1^r$ and summing over
$r\in\hf+\Z$, we obtain by using $\A(t)=-\sum_{r\in\hf+\Z}
t^r\gamma_r^-\gamma_{-r}^+$ (which follows from the definition of
$\A(t)$ and commutation relations among $\gamma_r^\pm$) that
$$
\begin{aligned}
\trace q^{L_0} & \A (qt_1)\A(t_2) \cdots \A(t_n) \\ &=
\sum_{s=0}^{n-1} (-1)^{s+1} \sum_{1<i_1<\cdots<i_s\leq n}
\mathfrak A^{(0)}_{-1}(q; t_1t_{i_1}\cdots
t_{i_s},\dots,\hat{t}_{i_1},\dots,\hat{t}_{i_s},\dots).
\end{aligned}$$
This finishes the proof.
\end{proof}
\subsection{The $q$-dimension formula of an $\hgl$-module of level $-l$}
Let $l\in\N$. We shall denote the $q$-dimension of the module
$L(\hgl; \Lambda(\la),-l)$ by
$$\mathsf Q_{-l}^\la(q) := \trace_{L(\hgl;\Lambda(\la),-l)}
q^{L_0}.
$$

\begin{lem} \label{gl:qdim:identity}
We have the following identity:
$$
\prod_{i=1}^l \frac{1}{(z_iq^\hf)_\infty(z_i^{-1}q^\hf)_\infty} =
\sum_{\la} \mathsf{ch}_\la^{gl}(z_1,\dots,z_l)\cdot
\trace_{L(\hgl;\Lambda(\la),-l)} q^{L_0}.
$$
\end{lem}

\begin{proof} Apply $\trace_{\F^{-l}}q^{L_0}z_1^{e_{11}}\cdots
z_l^{e_{ll}}$ to both sides of Proposition~\ref{gl:duality}.  From
the formulas for $e_{ii}$ in \cite{W1}, on the left-hand side,
$z_i^{e_{ii}}$ only acts on the $i^{th}$ tensor factor and the
resulting formula follows from the structure of the bosonic Fock
space $\F^{-1}$.  For the right-hand side, only
$z_1^{e_{11}}\cdots z_l^{e_{ll}}$ acts on $V_\la(GL(l))$.
\end{proof}

\begin{prop} \label{gl:qdim:1}
Let $k \in \Z$. The $q$-dimension of the irreducible $\hgl$-module
of highest weight $\Lambda(k)$ and level $-1$ is
$$
\mathsf Q_{-1}^{(k)}(q) = \frac{1}{(q)_\infty^2}\sum_{m\geq 0}
(-1)^m q^{\hf m(m+1)+ |k| (m+\half)}.
$$
\end{prop}

\begin{proof}
By Lemma \ref{gl:qdim:identity} and the identity in Theorem
\ref{identity:ff} with $u = zq^\hf$, we have
\begin{align*}
\sum_{k\in \Z} z^k \mathsf Q^{(k)}_{-1}(q) &=
\frac{1}{(zq^\hf)_\infty (z^{-1}q^\hf)_\infty} \\
 %
 &= \frac1{(q)_{\infty}^2}
\sum_{m=0}^\infty (-1)^m q^{\hf m(m+1)}
   ( \sum_{k\ge 0}   q^{k(m+\hf)} z^k
   + \sum_{k >0} q^{k(m+\hf)} z^{-k})\\
 &= \frac1{(q)_{\infty}^2}  \sum_{k\in \Z}
 \sum_{m=0}^\infty  z^{k} (-1)^m q^{\hf m(m+1)} q^{|k|(m+\hf)}.
\end{align*}
The proposition follows by comparing the coefficients of $z^k$ on
both sides.
\end{proof}

\begin{thm} \label{gl:qdim:l}
The $q$-dimension of the $\hgl$-module of highest weight
$\Lambda(\la)$ and level $-l$ is
$$
\mathsf Q_{-l}^\la(q) = \sum_{\sigma\in S_l} (-1)^{\ell(\sigma)}
\mathsf Q_{-1}^{(k_1)}(q)\cdots \mathsf Q_{-1}^{(k_l)}(q),
$$
where $k_i = (\la+\rho-\sigma(\rho),\varepsilon_i)$.
\end{thm}

\begin{proof}
Note that
$$\begin{aligned}
\prod_{i=1}^l \frac{1}{(z_iq^\hf)_\infty(z_i^{-1}q^\hf)_\infty}
&= \prod_{i=1}^l\left(\sum_{k\in\Z} z_i^k\mathsf
Q_{-1}^{(k)}(q)\right) \\
& = \sum_{\mathbf{k}\in\Z^l}z_1^{k_1}\cdots z_l^{k_l}\mathsf
Q_{-1}^{(k_1)}(q)\cdots \mathsf Q_{-1}^{(k_l)}(q).
\end{aligned}$$

Using Lemma \ref{gl:qdim:identity}, the character formula
(\ref{weylchar}) and the Weyl denominator formula, we get
$$\begin{aligned}
\sum_{\la} |z_j^{\la_i+l-i}|\, \mathsf Q_{-l}^\la(q)
&= |z_j^{l-i}|\sum_{\mathbf{k}\in\Z^l}z_1^{k_1}\cdots
z_l^{k_l}\mathsf Q_{-1}^{(k_1)}(q)\cdots \mathsf Q_{-1}^{(k_l)}(q) \\
&= \sum_{\sigma\in S_l}(-1)^{\ell(\sigma)}\prod_{i=1}^l
z_i^{(\sigma(\rho),\varepsilon_i)}\sum_{\mathbf{k}\in\Z^l}z_1^{k_1}\cdots
z_l^{k_l}\mathsf Q_{-1}^{(k_1)}(q)\cdots \mathsf
Q_{-1}^{(k_l)}(q).
\end{aligned}$$
Comparing the coefficients of $\prod_{i=1}^l z_i^{\la_i+l-i}$ on
both sides gives the result.
\end{proof}

\section{The $\cinf$-correlation functions and $q$-dimension formulas}
\label{sec:cinf}
\subsection{Lie algebra $\cinf$}

Let
\begin{eqnarray*}
 {\overline{c}}_{\infty}
      = \{ (a_{ij})_{i,j \in \Z} \in \gl\mid
                    a_{ij} = - (-1)^{i+j}a_{1-j,1-i} \}
\end{eqnarray*}
be a Lie subalgebra of $\gl$ of type $C$ \cite{DJKM2}. Denote by
$\cinf$ the central extension of ${\overline{c}}_{\infty}$ given by
the restriction of the $2$-cocycle (\ref{eq_cocy}) to
$\overline{c}_{\infty}$. Then $\cinf$ inherits from $\hgl$ a
triangular decomposition with Cartan subalgebra $(\cinf)_0$. We let
\begin{eqnarray*}
  H^c_i & = & E_{ii} + E_{-i, -i} - E_{i+1, i+1} - E_{1-i, 1-i},
   \quad i \in \N,    \\
   H^c_0 & = & E_{0,0} -E_{1,1} + C.
\end{eqnarray*}
Denote by $\hL^c_i \in (\cinf)_0^* $ the $i$-th fundamental weight
of $\cinf$, i.e. $\hL^c_i (H^c_j ) = \delta_{ij}$.

The Dynkin diagram of $\cinf$ is:
  \begin{equation*}
  \begin{picture}(150,45) 
  \put(10,20){$\circ$}
  \put(17,24){\line(1,0){32}}
  \put(17,22){\line(1,0){32}}
  \put(29,20){$>$}
  \put(50,20){$\circ$}
  \put(57,23){\line(1,0){32}}
  \put(91,20){$\circ$}
  \put(98,23){\line(1,0){32}}
  \put(131,20){$\circ$}
  \put(140,23){\dots}
  \put(10,9){$0$}
  \put(50,9){$1$}
  \put(91,9){$2$}
  \put(131,9){$3$}
  \end{picture}
  \end{equation*}

The Lie algebra $\cinf$ is generated by
\begin{equation*}
E_{ij}-(-1)^{i+j}E_{1-j,1-i},\quad i,j\in\Z,
\end{equation*}
which can be organized into a generating series as
$$
E(z,w):=\sum_{i,j}\big{(}E_{ij}-(-1)^{i+j}E_{1-j,1-i}\big{)}
z^{i-1}w^{-j}.
$$
Following \cite{TW}, we introduce the following operators in
$\cinf$:
$$
\begin{aligned}
\no\Ct(t)\no
&=\sum_{r\in\hf+\Z_+}(t^r-t^{-r}) \left(E_{r+\hf,r+\hf}-E_{\hf-r,\hf-r} \right),\\
\Ct(t) &=\no\Ct(t)\no+\frac{2}{t^\hf-t^{-\hf}}C.
\end{aligned}
$$
\subsection{The 1-point $\cinf$-function of level $-\hf$}
\label{section:onept:cinf}

Let $\chi(z)=\sum_{r\in\hf+\Z}\chi_r z^{-r-\hf}$ be a free bosonic
field with commutation relations
\begin{equation*}
[\chi_r,\chi_s]=(-1)^{r+\hf}\delta_{r+s,0},\quad r,s\in\hf+\Z.
\end{equation*}
The Fock space $\F^{-\hf}$ (cf. \cite{FF, W1}) generated by
$\chi(z)$ is a $\cinf$-module of level $-\hf$ given by
$$
E(z,w)=\no\chi(z)\chi(-w)\no,
$$
or equivalently by letting
$$
E_{ij}-(-1)^{i+j}E_{1-j,1-i}=(-1)^{-j}\no\chi_{-i+\hf}\chi_{j-\hf}\no.
$$
When acting on $\F^{-\hf}$ we have
\begin{align*}
\Ct(t) &=\sum_{r\in\hf+\Z_+}\left((-1)^{r+\hf}t^r\chi_{-r}\chi_r
-(-1)^{r+\hf}t^{-r}\chi_r\chi_{-r}\right)  \\
 &=\sum_{r\in\hf+\Z} (-1)^{r+\hf}t^r\chi_{-r}\chi_r,  \\
\no\Ct(t)\no&=\sum_{r\in\hf+\Z_+}
(-1)^{r+\hf}(t^r-t^{-r})\no\chi_{-r}\chi_r\no \\
&=\sum_{r\in\hf+\Z_+}\left((-1)^{r+\hf}t^r\chi_{-r}\chi_r
-(-1)^{r+\hf}t^{-r}\chi_{-r}\chi_r\right).
\end{align*}

\begin{lem}\label{333}
We have $[\Ct(t),\chi_s] = -(t^s-t^{-s}) \chi_s.$
\end{lem}

\begin{proof}
For $s>0$ we compute that
\begin{align*}
[\Ct(t),\chi_s] &= \sum_{r>0} (-1)^{r+\hf}t^r[\chi_{-r},\chi_s]
\chi_r -
(-1)^{r+\hf}t^{-r} \chi_r [\chi_{-r},\chi_s] \\
&  = -(t^s-t^{-s}) \chi_s.
\end{align*}
The case for $s<0$ is similar.
\end{proof}

Let $L_0$ be the zero mode of the Virasoro field so that
$[L_0,\chi_{-r}]=r\chi_{-r}$ and $L_0\vac = 0$. Define the
$n$-point $\cinf$-function of level $-\hf$ as
$$
\mathfrak C^{(0)}_{-\hf}(q; t_1,\dots,t_n) =
\trace_{\F^{-\hf}}q^{L_0}\Ct(t_1)\cdots \Ct(t_n).
$$

\begin{thm}  \label{th:C1ptHalf}
The $1$-point $\cinf$-function $\mathfrak{C}^{(0)}_{-\hf}(q;t)$ is
equal to
\begin{align*}
&\frac{1}{(q^\hf)_\infty(t^{-\hf} -t^\hf)} - \frac{1}{(q^\hf)_\infty
 (1-q^{-\hf})} \frac{t^\hf(tq^{\frac{3}{2}})_\infty}{(tq)_\infty}
{}_2\Phi_2(q^\hf,q^\hf;q^{\frac{3}{2}},tq^{\frac{3}{2}};q^2t)\\
&+\frac{1}{(q^\hf)_\infty
 (1 -q^{-\hf})}
\frac{t^{-\hf}(t^{-1}q^{\frac{3}{2}})_\infty}{(t^{-1}q)_\infty}
{}_2\Phi_2(q^\hf,q^\hf;q^{\frac{3}{2}},t^{-1}q^{\frac{3}{2}};q^2t^{-1}).
\end{align*}
\end{thm}

\begin{proof}
The Fock space $\F^{-\hf}$ has a basis given by
\begin{equation*}
v_\la =\chi_{-\la_1+\hf} \chi_{-\la_2+\hf} \cdots \vac,
\end{equation*}
where $\la=(\la_1,\la_2,\cdots)$ runs over all partitions. By Lemma
\ref{333}
\begin{equation*}
\Ct(t)v_\la =\left (
\sum_{i=1}^l(t^{\la_i-\hf}-t^{-\la_i+\hf})+\frac{1}{t^{-\hf}
-t^\hf}\right) v_\la,
\end{equation*}
and hence we have
\begin{equation*}
\mathfrak{C}^{(0)}_{-\hf}(q;t) =\frac{1}{(q^\hf)_\infty(t^{-\hf}
-t^\hf)}+\sum_{l=1}^\infty q^{-\frac{l}{2}}\sum_{l(\la)=l}
q^{|\la|}\sum_{i=1}^l(t^{\la_i-\hf}-t^{-\la_i+\hf}).
\end{equation*}
Now the theorem follows by applying (\ref{555}) (with $x=1$) twice.
\end{proof}

We remark that the above $q$-hypergeometric series is again of type
II.

\subsection{A $q$-difference equation for $\cinf$-functions of level $-\hf$}

\begin{thm}
The $n$-point $\cinf$-function satisfies the $q$-difference
equation:

$\mathfrak C^{(0)}_{-\hf}(q; qt_1,t_2,\dots,t_n) = $
$$\sum_{s=0}^{n-1}\sum_{1<i_1<\cdots<i_s\leq
n}\sum_{\epsilon_{i_a}=\pm 1}(-1)^{s+\#\epsilon} \mathfrak
C^{(0)}_{-\hf}(q; t_1t_{i_1}^{\epsilon_{i_1}}\cdots
t_{i_s}^{\epsilon_{i_s}},\dots,\widehat{t}_{i_1},\dots,\widehat{t}_{i_s},\dots),$$
where $\#\epsilon$ stands for the number of $\epsilon_{i_a}$'s that
are equal to $-1$.
\end{thm}

\begin{proof}
For $r > 0$, by Lemma \ref{333} we have
\begin{equation}\label{cinf:qdiff:one}
\Ct(t_2)\cdots \Ct(t_n)\chi_r =
\sum_{S\subset\{2,\dots,n\}}(-1)^{|S|}\prod_{i\in
S}(t_i^r-t_i^{-r})\chi_r\prod_{i\not\in S}\Ct(t_i),
\end{equation}
and
\begin{equation} \label{cinf:qdiff:two}
\Ct(t_2)\cdots \Ct(t_n)\chi_{-r} =
\sum_{S\subset\{2,\dots,n\}}(-1)^{|S|}\prod_{i\in
S}(t_i^{-r}-t_i^{r})\chi_{-r}\prod_{i\not\in
S}\Ct(t_i).\end{equation}

We will write $\trace$ instead of $\trace_{\F^{-\hf}}$ below. Let
us simplify notation by setting $$S_r^\pm := \trace
q^{L_0}\chi_{\mp r} \Ct(t_2)\cdots \Ct(t_n)\chi_{\pm r}.$$ Then
(\ref{cinf:qdiff:one}) and (\ref{cinf:qdiff:two}) imply the
following:

\begin{equation} \label{cinf:qdiff:trone} S_r^+ =
\sum_{S\subset\{2,\dots,n\}}(-1)^{|S|}\prod_{i\in
S}(t_i^r-t_i^{-r})\trace q^{L_0}\chi_{-r}\chi_r\prod_{i\not\in
S}\Ct(t_i),\end{equation} and \begin{equation}
\label{cinf:qdiff:trtwo} S_r^- =
\sum_{S\subset\{2,\dots,n\}}(-1)^{|S|}\prod_{i\in
S}(t_i^{-r}-t_i^{r})\trace q^{L_0}\chi_{r}\chi_{-r}\prod_{i\not\in
S}\Ct(t_i).\end{equation}

It is clear that
$$\prod_{i\in
S}(t_i^{\pm r}-t_i^{\mp r}) = \sum_{\epsilon_i=\pm 1, i\in
S}(-1)^{\#\epsilon}\prod_{i\in S}(t_i^{\epsilon_{i}})^{\pm r}.$$

We will compute $$\sum_{r>0} \left((-1)^{r+\hf}t_1^rS_r^+ -
(-1)^{r+\hf} t_1^{-r}S_r^-\right)$$ in two different ways.  Using
(\ref{cinf:qdiff:trone}) and (\ref{cinf:qdiff:trtwo}) it is equal to

$$\begin{aligned}
\sum_{r>0}&\sum_{S\subset\{2,\dots,n\}}(-1)^{|S|}\sum_{\epsilon_{i}=\pm
1,i\in S}(-1)^{\#\epsilon}\trace q^{L_0}\\ &\times
\left((-1)^{r+\hf}(t_1\prod_{i\in
S}(t_i^{\epsilon_{i}}))^r\chi_{-r}\chi_r-(-1)^{r+\hf}(t_1\prod_{i\in
S}(t_i^{\epsilon_{i}}))^{-r}\chi_{r}\chi_{-r}\right)\prod_{i\not\in
S}\Ct(t_i),\end{aligned}$$ which, using the definition of $\Ct(t)$,
is
\begin{equation} \label{cinf:qdiff:rhs}
\sum_{s=0}^{n-1}(-1)^s\sum_{1<i_1<\cdots<i_s\leq
n}\sum_{\epsilon_{i_a}=\pm 1}(-1)^{\#\epsilon} \mathfrak
C^{(0)}_{-\hf}(t_1t_{i_1}^{\epsilon_{i_1}}\cdots
t_{i_s}^{\epsilon_{i_s}},\dots,\widehat{t}_{i_1},\dots,\widehat{t}_{i_s},\dots).
\end{equation}

On the other hand, using the commutation relation $[L_0,\chi_{-r}]
= r\chi_{-r}$, we see that for any $r$,
$$
\begin{aligned}
\trace q^{L_0}\chi_{-r} \Ct(t_2)\cdots \Ct(t_n)\chi_r  &=
\trace \chi_r q^{L_0}\chi_{-r}\Ct(t_2)\cdots \Ct(t_n) \\
&= q^r\trace q^{L_0}\chi_r\chi_{-r} \Ct(t_2)\cdots \Ct(t_n).
\end{aligned}
$$
Now we have that
\begin{equation}\label{cinf:qdiff:lhs}
\sum_{r>0} \left((-1)^{r+\hf}t_1^rS_r^+ -
(-1)^{r+\hf}t_1^{-r}S_r^-\right) = \trace q^{L_0}
\tilde{\Ct}(qt)\Ct(t_2)\cdots \Ct(t_n),
\end{equation}
where $$\tilde{\Ct}(t)=\sum_{r>0}
\left((-1)^{r+\hf}t^r\chi_r\chi_{-r} -
(-1)^{r+\hf}t^{-r}\chi_{-r}\chi_r\right).
$$
One verifies that
$\tilde{\Ct}(t)=\Ct(t)$, by observing that $\sum_{r>0} t^r
=\frac{t^{\hf}}{1-t}$ coincides with  $- \sum_{r>0} t^{-r}
=\frac{t^{-\hf}}{t^{-1} -1}$. Equating (\ref{cinf:qdiff:rhs}) and
(\ref{cinf:qdiff:lhs}) now completes the proof.
\end{proof}
\subsection{The $n$-point $\cinf$-functions of level $l-\hf$}

Let $l \in \N$. Let $\F^l$ be the Fock space generated by $l$ pairs
of free fermions $\psi^{\pm,p}(z)=\sum_{r\in\hf+\Z}\psi^{\pm,p}_r
z^{-r-\hf}$, $p=1,\cdots,l$, with non-trivial commutation relations
$$
[\psi^{+,p}_r,\psi^{-,q}_s]_+=\delta_{p,q}\delta_{r+s,0}, \quad \text{ for }
r,s\in\hf+\Z.
$$
Let $\F^{l-\hf} \equiv \F^{-\hf}\otimes\F^l$. The Lie algebra
$\cinf$ acts on $\F^{l-\hf}$ by (see \cite{W1})
$$\begin{aligned}
\sum_{i,j\in\Z} &(E_{i,j} - (-1)^{i+j}E_{1-j,1-i})z^{i-1}w^{-j} \\
&= \sum_{k=1}^l \left(\no\psi^{+,k}(z)\psi^{-,k}(w)\no +
\no\psi^{+,k}(-w)\psi^{-,k}(-z)\no\right) + \no\chi(z)\chi(-w)\no.
\end{aligned}
$$
It follows that the operator $\Ct(t)$ acts on $\F^{l-\hf}$ as
$$
\Ct(t) = \sum_{r \in\hf+\Z} \sum_{p=1}^l t^r \left(
\psi_{-r}^{+,p}\psi_r^{-,p}+\psi_{-r}^{-,p}\psi_r^{+,p}\right) +
\sum_{r\in\hf+\Z} (-1)^{r+\hf} t^r\chi_{-r}\chi_r.$$

{F}rom \cite[Lemma 4.6]{W1} we have an action of the Lie
superalgebra $\mathfrak{osp}(1,2l)$ on $\F^{l-\hf}$ that commutes
with the action of $\cinf$. It is known that the finite-dimensional
simple $\mathfrak{osp}(1,2l)$-modules, denoted by
$V_\la(\mathfrak{osp}(1,2l))$, are parameterized by partitions $\la$
of length $\leq l$ as highest weights. Associated to a partition
$\la = (\la_1, \ldots, \la_l)$ we define a highest weight
$\Lambda(\la)$ of $\cinf$ via
$$\Lambda (\lambda) =
   (l - \hf -j) \hL^c_0 + \sum_{k =1}^j \hL^c_{\la_k},$$
if $\la_1 \geq \ldots \geq \la_j > \la_{j+1} = \ldots = \la_l =0.$

\begin{prop}\label{444} \cite[Theorem 4.3]{W1}
We have the $(\mathfrak{osp}(1,2l),\cinf)$-module decomposition
$$
\F^{l-\hf} \cong \bigoplus_{\la}
V_\la(\mathfrak{osp}(1,2l))\otimes L(\cinf;\Lambda(\la), l-\hf).
$$
\end{prop}

Let ${\bf t}=(t_1,\ldots,t_n)$. The Bloch-Okounkov $n$-point
$\cinf$-function of level $l-\hf$ (associated to a partition $\la$
of length $\leq l$) is defined as
$$
\mathfrak C^\la_{l-\hf}(q; {\bf t}) \equiv \mathfrak
C^\la_{l-\hf}(q; t_1,\dots,t_n) =
\trace_{L(\cinf;\Lambda(\la),-l)}q^{L_0}\Ct (t_1) \Ct (t_2)\cdots\Ct
(t_n).
$$

\begin{lem} \label{cinf:charequiv}
The character $\mathsf{ch}_\la^{\mathfrak{osp}}(z_1,\dots,z_l)
:=\trace_{V_\la(\mathfrak{osp}(1,2l))} (z_1^{e_{11}}\cdots
z_l^{e_{ll}})$ of the irreducible $\mathfrak{osp}(1,2l)$-module of
highest weight $\la$ is given by
$$
\mathsf{ch}_\la^{\mathfrak{osp}}(z_1,\dots,z_l) = \frac{\left
|z_j^{\la_i+l-i+\hf} -z_j^{-(\la_i+l-i+\hf)} \right |}{\left
|z_j^{l-i+\hf}-z_j^{-(l-i+\hf)} \right |}.
$$
\end{lem}

\begin{proof}
Recall (cf.~\cite{K2}) that the positive even roots for
$\mathfrak{osp}(1,2l)$ are given by $\Phi_0 =
\{2\varepsilon_i,\varepsilon_i\pm \varepsilon_j |1\le i\not= j \le
l\}$ and the positive odd roots by $\Phi_1 = \{\varepsilon_i |
1\le i \le l\}$. The character formula for
$V_\la(\mathfrak{osp}(1,2l))$ is given by \cite{K2}
$$
\mathsf{ch} V_\la(\mathfrak{osp}(1,2l)) = \frac{\displaystyle
\prod_{\alpha\in\Phi_1} (1+e^{-\alpha})\displaystyle
\sum_{\sigma\in W(C_l)}(-1)^{\ell(\sigma)}
e^{\sigma(\la+\rho_{osp})-\rho_{osp}}}{\displaystyle\prod_{\alpha\in\Phi_0}(1-e^{-\alpha})},
$$
where $\rho_{osp} = \hf\sum_{\alpha\in\Phi_0} \alpha -
\hf\sum_{\alpha\in\Phi_1} \alpha$.

Denote by $\Phi_B$ the root system for type $B$, i.e.~$\Phi_B =
\{\varepsilon_i,\varepsilon_i\pm \varepsilon_j | i\not= j\}$. It is
easy to see that $\rho_{osp}$ is equal to $\rho_B$.  Note that
$$
\frac{1+e^{-\varepsilon_i}}{1-e^{-2\varepsilon_i}} =
\frac{1}{1-e^{-\varepsilon_i}}.
$$
Since the Weyl group is the same for type $C$ and type $B$, the
character of $V_\la(\mathfrak{osp}(1,2l))$ coincides with
irreducible $\mathfrak{so}(2l+1)$-character of highest weight
$\la$, which is known (cf. \cite{FH}) to be given by the
right-hand side of the formula in the lemma.
\end{proof}

Let use denote by $\mathbf{F}(z,q,\mathbf{t}) = \trace_{\F^{1}}
z^{e_{11}}q^{L_0}\Ct(t_1)\cdots \Ct(t_n)$.

\begin{lem} \label{cinf:qseries:one} We have the following $q$-series identity:
$$\begin{aligned} \trace_{\F^{-\hf}}q^{L_0}\Ct(t_1)\cdots \Ct(t_n)\cdot \prod_{i=1}^l
\mathbf{F}(z_i,q,\mathbf{t}) = \sum_{\la}
\mathsf{ch}_\la^{osp}(z_1,\dots,z_l) \mathfrak C^\la_{l-\hf}(q;
\mathbf{t}).
\end{aligned}$$
\end{lem}

\begin{proof}
The identity results by applying $\trace_{\F^{l-\hf}}
z_1^{e_{11}}\cdots z_l^{e_{ll}}q^{L_0} \Ct(t_1)\cdots \Ct(t_n)$ to
both sides of Proposition~\ref{444}. On the left-hand side,
$z_i^{e_{ii}}$ only acts on the $i^{th}$ tensor factor of $\F^l$
and not on $\F^{-\hf}$. For the right-hand side,
$z_1^{e_{11}}\cdots z_l^{e_{ll}}$ acts only on the first tensor
factor of $V_\la(\mathfrak{osp}(1,2l))$ and $q^{L_0}\Ct(t_1)\cdots
\Ct(t_n)$ acts only on the second tensor factor.\end{proof}

Define the theta function
\begin{eqnarray*}
\Theta (t) :=
 (t^{\hf} -t^{-\hf})(q)_\infty^{-2} (qt)_\infty(qt^{-1})_\infty,
\end{eqnarray*}
and let $\Theta^{(k)} (t) = \left(t\frac{d}{dt} \right)^k \Theta
(t)$ for $k\in\Z_+.$ Denote
\begin{align*}
 F_{bo}(q;{\bf t}) := \frac{1}{(q)_\infty}
\sum_{\sigma\in S_n} \frac{{\rm det}
\Big{(}\frac{\Theta^{(j-i+1)}(t_{\sigma(1)}\cdots
t_{\sigma(n-j)})}{(j-i+1)!} \Big{)}_{i,j=1}^n}
{\Theta(t_{\sigma(1)}) \Theta(t_{\sigma(1)}t_{\sigma(2)})\cdots
\Theta(t_{\sigma(1)}t_{\sigma(2)}\cdots t_{\sigma(n)})}
\end{align*}
where it is understood below that $\frac{1}{(-k)!}=0$, for
$k\in\N$.
We recall the formula for the original Bloch-Okounkov $n$-point
correlation function of the $\ainf$-module $\F^1_{(k)}$ of level
$1$ \cite{BO, Ok} is given by
\begin{equation}\label{zzz}
\trace_{\F^{1}_{(k)}}(q^{L_0}\A(t_1) \cdots \A(t_n)) =
q^{\frac{k^2}{2}}(t_1\cdots t_n)^k F_{bo}(q;t_1,\cdots,t_n).
\end{equation}


\begin{thm}  \label{th:Cfunction}
The $n$-point $\cinf$-function, $\mathfrak{C}^\la_{l-\hf}(q;
t_1,\dots,t_n)$, is equal to
\begin{eqnarray*}
\mathfrak{C}^{(0)}_{-\hf}(q;\mathbf{t})\cdot\sum_{\sigma\in
W(B_l)}\left(-1\right)^{\ell(\sigma)} q^{\frac{\Vert
\la+\rho_B-\sigma(\rho_B) \Vert^2}{2}}
 \prod_{a=1}^l
 \Big(
\sum_{\vec{\epsilon}_a\in\{\pm 1\}^n}
 [\vec{\epsilon}_a] (\Pi {\bf t}^{\vec{\epsilon}_a})^{k_a}
 F_{bo}(q;{\bf t}^{{\vec{\epsilon}_a}})
 \Big),
\end{eqnarray*}
where $k_a=(\la+\rho_B-\sigma(\rho_B),\varepsilon_a)$, and for an
expression $\vec{\epsilon}=(\epsilon_1,\cdots,\epsilon_n)$ we define
$[\vec{\epsilon}]=\epsilon_1\cdots\epsilon_n$ and $\prod
\mathbf{t}^{\vec{\epsilon_a}} = t_1^{\epsilon_1}\cdots
t_n^{\epsilon_n}$.
\end{thm}

\begin{proof}
This proof mirrors the one used in \cite{TW}. The Weyl denominator
of type $B_l$ (also the denominator of $\mathsf{ch}_\la^{osp}$)
reads that
$$
|z_j^{l-i+\hf}-z_j^{-(l-i+\hf)}| = \sum_{\sigma\in W(B_l)}
(-1)^{\ell(\sigma)}\mathbf{z}^{\sigma(\rho_B)},
$$
where we denote $\mathbf{z}^\mu = z_1^{\mu_1}\cdots z_l^{\mu_l}$,
for $\mu = (\mu_1,\dots,\mu_l)$. It follows from Lemma
\ref{cinf:qseries:one} and $\Ct(t)=\A(t)-\A(t^{-1})$ that
$$
\begin{aligned}
\sum_{\sigma\in W(B_l)}
&(-1)^{\ell(\sigma)}\mathbf{z}^{\sigma(\rho_B)}\cdot
\mathfrak{C}^{(0)}_{-\hf}(q; \mathbf{t}) \\
&\qquad\qquad\times \prod_{a=1}^l
 \left( \sum_{k_a\in\Z} z_a^{k_a} q^{\frac{k_a^2}{2}}
 \sum_{\vec{\epsilon}_a  \in \{\pm 1\}^n}  [\vec{\epsilon}_a] \cdot
 \left(\Pi{\bf t}^{\vec{\epsilon}_a} \right)^{k_a}
 F_{bo}(q;{\bf t}^{{\vec{\epsilon}_a}})
 \right) \\
&= \sum_{\la} |z_j^{\la_i+l-i+\hf}-z_j^{-(\la_i+l-i+\hf)}|\cdot
 \mathfrak{C}^\la_{l-\hf}(q; \mathbf{t}),
 \end{aligned}
 $$
where we have used (\ref{zzz}) in the above calculation. Among the
monomials $\mathbf{z}^\mu$ in the expansion of
$|z_j^{\la_i+l-i+\hf}-z_j^{-(\la_i+l-i+\hf)}|$, there is exactly one
dominant monomial with $\mu_1\geq\cdots\geq\mu_l\geq 0$, which is
$\mathbf{z}^{\la+\rho_B}$.  The theorem follows by comparing the
coefficient of $\mathbf{z}^{\la+\rho_B}$ on both sides of the above
equation.
\end{proof}
\subsection{The $q$-dimension of a $\cinf$-module of level
$l-\half$}

Let $l \in \N$. The $q$-dimension of the $\cinf$-module
$L(\cinf;\Lambda(\la),l-\hf)$ is
$${}^c\mathsf{Q}_{l-\hf}^\la(q) :=
\trace_{L(\cinf;\Lambda(\la),l-\hf)} q^{L_0}.
$$

We can derive the following $q$-dimension formula from the
$(\mathfrak{osp}(1,2l), \cinf)$-duality (see
Proposition~\ref{444}), where the second formula is easily seen to
be equivalent to the first by \cite[Lemma 3.10]{TW}.

\begin{thm}
Let $\la =(\la_1,\ldots, \la_l)$ be a partition of length $\le l$.
We have
\begin{align*}
{}^c\mathsf{Q}_{l-\hf}^\la(q)
 & =
\frac{1}{(q^{\hf})_\infty (q)_\infty^l} \cdot
\sum_{\sigma\in W(B_l)}(-1)^{\ell(\sigma)} q^{\frac{\Vert
\la+\rho_B-\sigma(\rho_B) \Vert2}{2}} \\
 &=
\frac{1}{(q^{\hf})_\infty (q)_\infty^l} \cdot
q^{\frac{\Vert\la\Vert2}{2}}
  \prod_{1\leq i\leq l} \left(1-q^{\la_i+l-i+\hf}\right) \times\\
 &\quad \times \prod_{1\leq i<j\leq l}
\left(1-q^{\la_i-\la_j+j-i}\right)
\left(1-q^{\la_i+\la_j+2l-i-j+1}\right).
\end{align*}
\end{thm}

\begin{proof}
By applying $\trace z_1^{e_{11}}\cdots z_l^{e_{ll}} q^{L_0}$ to
both sides of the duality in Proposition~\ref{444}, we obtain
$$\prod_{i=1}^l \left(\trace_{\F^1} z_i^{e_{ii}}q^{L_0}\right)
\trace_{\F^{-\half}} q^{L_0} = \sum_{\la}
\mathsf{ch}_\la^{osp}(z_1,\dots,z_l)
{}^c\mathsf{Q}_\la^{l-\half}(q).
$$
Noting by the Jacobi triple product identity that
$$
\trace_{\F^1} z^{e_{ii}}q^{L_0}
 = \prod_{r \ge 0} (1 +q^{r+\hf} z) (1 +q^{r+\hf} z^{-1})
  = \frac{1}{(q)_\infty}\sum_{k\in\Z}
z^kq^{\frac{k^2}{2}}
$$
and that
$$
\trace_{\F^{-\half}} q^{L_0} = \frac{1}{(q^\half)_\infty},
$$
a completely analogous argument as for Theorem~\ref{th:Cfunction}
applies.
\end{proof}

\subsection{The $n$-point $\cinf$-functions of level $-l$}

Let $l \in \N$. Let $\F^{-l}$ denote the Fock space generated $l$
pairs of free bosonic fields
$\gamma^{\pm,p}(z)=\sum_{r\in\hf+\Z}\gamma^{\pm,p}_r z^{-r-\hf}$
($p=1,\ldots,l$) with non-trivial commutation relations
$$[\gamma^{+,p}_r,\gamma^{-,q}_s]=\delta_{pq}\delta_{r+s,0}, \quad \text{ for }
r,s\in\hf+\Z.
$$
According to \cite[Section 5.2]{W1} there is an action of $\cinf$
on $\F^{-l}$, from which we conclude that
$$\Ct(t) =  \sum_{r\in\Z+\half} \sum_{p=1}^l
t^r (- \gamma_{-r}^{+,p}\gamma_r^{-,p} +
\gamma_{-r}^{-,p}\gamma_r^{+,p}).
$$
There is also an action of the Lie group $O(2l)$ on $\F^{-l}$ which
commutes with the action of $\cinf$.

Denote the parameter set for simple $O(2l)$-modules (cf.
\cite{BtD, W1}) by
\begin{eqnarray*}
  \mathcal C :=
   & & \left\{ (\la_1, \la_2, \ldots, \la_l ) \mid
       \la_1 \geq \la_2 \geq \cdots \geq \la_l > 0, \la_i \in \Z   \right \}
                               \\
   & & \cup   \left\{ (\la_1, \la_2, \ldots, \la_{l-1}, 0 ) \otimes \det,
   (\la_1, \la_2, \ldots, \la_{l-1}, 0 ) \mid \right.\\
   & &   \left. \quad  \la_1 \geq \cdots
   \geq  \la_{l-1} \geq 0, \la_i \in \Z \right\}.
\end{eqnarray*}

Associated to $\la\in\mathcal C$ we define a highest weight
$\Lambda(\la)$ for $\cinf$
$$
\Lambda(\la) = (-l-\la_1)\Lambda_0^c + \sum_{k=1}^l
(\la_k-\la_{k+1})\Lambda_0^c,
$$
if $\la = (\la_1,\dots,\la_l)$ with $\la_l
> 0$, $$\Lambda(\la) = (-l-\la_1)\Lambda_0^c + \sum_{k=1}^j (\la_k -
\la_{k+1})\Lambda_0^c,$$ if $\la = (\la_1,\dots,\la_j,0,\dots,0)$,
and
$$
\Lambda(\la) = (-l-\la_1)\Lambda_0^c + \sum_{k=1}^{j-1}
(\la_k-\la_{k+1})\Lambda_0^c + (\la_j-1)\Lambda_j^c +
\Lambda_{2l-j}^c,
$$
if $\la = (\la_1,\dots,\la_j,0,\dots,0)\otimes\det$.

\begin{prop} \label{cinf:integral:duality} \cite[Theorem 5.3]{W1}
We have the following $(O(2l),\cinf)$-module decomposition:
$$
\F^{-l} \cong \bigoplus_{\la\in\mathcal C} V_\la(O(2l))\otimes
L(\cinf;\Lambda(\la),-l),
$$
where $V_\la(O(2l))$ is the irreducible $O(2l)$-module
parameterized by $\la$.
\end{prop}

\begin{definition}
The $n$-point $\cinf$-function of level $-l$ (associated to a
partition $\la$ of length $\leq l$) is
$$
\mathfrak C_{-l}^\la(q;t_1,\dots,t_n) =
\begin{cases}
\trace_{L(\cinf;\Lambda(\la),-l)}
q^{L_0}\Ct(t_1)\cdots\Ct(t_n), & \la_l \not= 0, \\
\trace_{L(\cinf;\Lambda(\la),-l)\oplus
L(\cinf;\Lambda(\la\otimes\det),-l)} q^{L_0}\Ct(t_1)\cdots\Ct(t_n),
& \la_l = 0.
\end{cases}
$$
\end{definition}

\begin{prop} \label{cinf:ctoa}
The $n$-point function of level $-1$, $\mathfrak
C^{(m)}_{-1}(q;t_1,\dots,t_n)$, is given by
$$[z^m] \sum_{\vec{\epsilon}_a\in\{\pm 1\}^n}
[\vec{\epsilon}_a]\trace_{\F^{-1}}z^{e_{11}}q^{L_0}\A(t_1^{\epsilon_1})\cdots
\A(t_n^{\epsilon_n}).$$
\end{prop}

\begin{proof}
Since by definition $\Ct(t) = \A(t) - \A(t^{-1})$, we have
$$\begin{aligned} \Ct(t_1)\cdots \Ct(t_n)
 &= \prod_{j=1}^n \left(\A(t_j) -\A(t_j^{-1})\right) \\
& = \sum_{\vec{\eps}\in\{\pm 1\}^n} \eps_1\eps_2 \cdots \eps_n
\A(t_1^{\eps_1}) \A(t_2^{\eps_2})\cdots \A(t_n^{\eps_n}).
\end{aligned}
$$
Recall that
$$
\F^{-1} \cong \oplus_{m\in\Z}\F^{-1}_{(m)}
$$
where $\F^{-1}_{(m)}$ is the $m$-eigenspace of $e_{11}$.
Proposition~\ref{cinf:integral:duality} states that
$\F^{-1}_{(m)}\cong \F^{-1}_{(-m)}$ as $c_\infty$-modules for
$m\not=0$, and also that $\F^{-1}_{(0)}\cong
L(c_\infty;\La(\emptyset),-1)\oplus
L(c_\infty;\La(\emptyset\otimes{\rm det}),-1)$. Noting that
$$
\begin{aligned}
\trace_{\F^{-1}} z^{e_{11}}q^{L_0}
\Ct(t_1)\cdots\Ct(t_n) &= \sum_{m\in\Z} z^m \trace_{\F^{-1}_{(m)}}
q^{L_0} \Ct(t_1)\cdots\Ct(t_n)
\\ &= \sum_{m\in\Z} z^m \mathfrak
C^{(|m|)}_{-1}(q;t_1,\dots,t_n),
\end{aligned}
$$
the result follows.
\end{proof}

\begin{lem} \cite[Lemma 3.2]{TW} \label{lem:charD}
Denote by $\mathsf{ch}^o_\la (z_1,\dots,z_l)$ the character of the
irreducible $O(2l)$-module $V_\la(O(2l))$. Then
\begin{eqnarray*}
 \mathsf{ch}^o_\la (z_1,\dots,z_l) = c_\la
\frac{\left |z_j^{\la_i+l-i} + z_j^{-(\la_i+l-i)}\right
|}{\left|z_j^{l-i}+z_j^{-(l-i)}\right |},
\end{eqnarray*}
where $c_\la =1$ if $\la_l =0$, and $c_\la =2$ if $\la_l \not= 0.$
\end{lem}

\begin{thm} \label{cinf:thm:integral}
Let $\la = (\la_1,\dots,\la_l)$ be a partition. The function
$\mathfrak C^\la_{-l}(q;t_1,\dots,t_n)$ is equal to
$$
\sum_{\sigma\in W(D_l)} (-1)^{\ell(\sigma)}\mathfrak
C^{(k_1)}_{-1}(q;t_1,\dots,t_n)\cdots \mathfrak
C^{(k_l)}_{-1}(q;t_1,\dots,t_n),
$$
where $k_i = (\la+\rho-\sigma(\rho),\varepsilon_i) \ge 0$.
\end{thm}

\begin{proof} Apply the trace of $z_1^{e_{11}}\cdots
z_l^{e_{ll}}q^{L_0}\Ct(t_1)\cdots\Ct(t_n)$ to both sides of the
isomorphism in Proposition~\ref{cinf:integral:duality}.
On the left-hand side, $z_i^{e_{ii}}$ acts only on the $i^{th}$
tensor factor, so it becomes
$$
\begin{aligned}
\prod_{i=1}^l \trace_{\F^{-1}}z_i^{e_{ii}}q^{L_0}&\Ct(t_1)\cdots
\Ct(t_n) = \prod_{i=1}^l \sum_{m_i\in\Z} z_i^{m_i} \mathfrak
C_{-1}^{(|m_i|)}(q;t_1,\dots,t_n) \\ &= \sum_{\mathbf{m}\in\Z^l}
z_1^{m_1}\cdots z_l^{m_l}\mathfrak
C_{-1}^{(|m_1|)}(q;t_1,\dots,t_n)\cdots \mathfrak
C_{-1}^{(|m_l|)}(q;t_1,\dots,t_n).
\end{aligned}
$$
On the right-hand side, $z_1^{e_{11}}\cdots z_l^{e_{ll}}$ acts only
on the module $V_\la(O(2l))$ while $q^{L_0}\Ct(t_1)\cdots \Ct(t_n)$
acts on $L(\cinf;\Lambda(\la),-l)$, so it becomes
$$
\begin{aligned}
\sum_{\la\in\mathcal C} \mathsf{ch}_\la^o(z_1,\dots,z_l)& \,
\mathfrak C_{-l}^\la(q;t_1,\dots,t_n) \\ &= \sum_{\la\in\mathcal
C} c_\la  \frac{\left|z_j^{\la_i+l-i} +
z_j^{-(\la_i+l-i)}\right|}{\left |z_j^{l-i}+z_j^{-(l-i)}\right|}
\cdot \mathfrak C_{-l}^\la(q,t_1,\dots,t_n).
\end{aligned}
$$

Thus, equating both sides gives us
$$
\begin{aligned}
\sum_{\mathbf{m}\in\Z^l} & z_1^{m_1}\cdots z_l^{m_l}\mathfrak
C_{-1}^{(|m_1|)}(q;t_1,\dots,t_n)\cdots \mathfrak
C_{-1}^{(|m_l|)}(q;t_1,\dots,t_n)\\
&= \sum_{\la\in\mathcal C} c_\la  \frac{\left|z_j^{\la_i+l-i} +
z_j^{-(\la_i+l-i)}\right|}{\left |z_j^{l-i}+z_j^{-(l-i)}\right|}
\cdot \mathfrak C_{-l}^\la(q,t_1,\dots,t_n).
\end{aligned}
$$
Multiplying both sides by the Weyl denominator of type $D_l$,
i.e.,
$$\frac12 |z_j^{l-i}+z_j^{-(l-i)}| = \sum_{\sigma\in W(D_l)}
(-1)^{\ell(\sigma)}\mathbf{z}^{\sigma(\rho)},
$$
we obtain that
$$
\begin{aligned}
\sum_{\sigma\in W(D_l)}
(-1)^{\ell(\sigma)}\mathbf{z}^{\sigma(\rho)} &
\sum_{\mathbf{m}\in\Z^l} z_1^{m_1} \cdots z_l^{m_l}\mathfrak
C_{-1}^{(|m_1|)}(q;t_1,\dots,t_n)\cdots \mathfrak
C_{-1}^{(|m_l|)}(q;t_1,\dots,t_n) \\
&= \sum_{\la\in\mathcal C} c_\la/2 \cdot \left|z_j^{\la_i+l-i} +
z_j^{-(\la_i+l-i)}\right| \mathfrak
C_{-l}^\la(q;t_1,\dots,t_n).
\end{aligned}
$$
The result follows by comparing the coefficients of $\prod_{i=1}^l
z_i^{\la_i+l-i}$ on both sides, noting that its coefficient in
$|z_j^{\la_i+l-i} + z_j^{-(\la_i+l-i)}|$ is precisely $2/c_\la$.
\end{proof}
The $1$-point (respectively, $2$-point) $\cinf$-function of level
$-l$ now follows by combining Proposition~\ref{cinf:ctoa},
Theorem~\ref{cinf:thm:integral}, Remark~\ref{rem:residue}
(respectively, Remark~\ref{rem:residue2}).
\subsection{The $q$-dimension of a $\cinf$-module of level $-l$}

Let us denote by ${}^c\mathsf Q_{-l}^\la(q)$ the $q$-dimension of
$L(\cinf;\Lambda(\la),-l)$, or in the case that $\la_l = 0$, the
$q$-dimension of $L(\cinf;\Lambda(\la),-l)\oplus L(\cinf;
\Lambda(\la \otimes\det),-l)$.

\begin{prop}
We have the following $q$-dimension formula of level $-1$ (for $k
\in \Z_+$):
$$
{}^c\mathsf Q_{-1}^{(k)}(q) = \frac{1}{(q)_\infty^2}\sum_{m\geq 0}
(-1)^m q^{\hf m(m+1)+ |k| (m+\half)}.
$$
\end{prop}

\begin{proof}
We can identify $L(\cinf;\Lambda^c(k),-1) =L(\ainf;\Lambda^a(k),-1)$
in $\F^{-1}$ for $k>0$ by comparing Propositions~\ref{gl:duality}
and \ref{cinf:integral:duality}, where we have temporarily used the
superscripts $a,c$ to distinguish the weights for $\ainf$ and
$\cinf$ respectively; also we have $L(\cinf;\Lambda^c(0),-1) \oplus
L(\cinf;\Lambda^c((0)\otimes\det),-1) =L(\ainf;\Lambda^a(0),-1)$.
Now the result follows from Proposition~\ref{gl:qdim:1}.
\end{proof}

\begin{thm}
Let $\la = (\la_1,\dots,\la_l)$ be a partition. We have the
following $q$-dimension formula:
$$
{}^c\mathsf Q_{-l}^\la(q) = \sum_{\sigma\in W(D_l)}
(-1)^{\ell(\sigma)}\, {}^c\mathsf Q_{-1}^{(k_1)}(q)\cdots
{}^c\mathsf Q_{-1}^{(k_l)}(q),
$$
where $k_i = (\la+\rho-\sigma(\rho),\varepsilon_i)$.
\end{thm}

\begin{proof}
The proof is similar to that of Theorem \ref{cinf:thm:integral}, and
we omit the details.
\end{proof}
\subsection{The $n$-point $\cinf$-functions of level $-l-\half$}

Denote by $\mathcal C_\half$ the following parameter set for
irreducible $O(2l+1)$-modules (\cite{BtD, W1}):
\begin{eqnarray*}
    \left\{ (\la_1, \la_2, \ldots, \la_l ) \otimes \det,
   (\la_1, \la_2, \ldots, \la_l ) \mid
    \la_1 \geq \la_2 \geq \ldots
   \geq  \la_l \geq 0, \la_i \in \Z \right\}.
\end{eqnarray*}

Associated to $\la \in \mathcal C_\half$  we define the highest
weight $\Lambda(\la)$ of $\cinf$  as
$$
\Lambda(\la) = (-l-\la_1-\half)\Lambda_0^c + \sum_{k=1}^j
(\la_k-\la_{k+1})\Lambda_k^c,
$$
if $\la = (\la_1,\dots,\la_j,0,\dots,0)$ and
$$
\Lambda(\la) = (-l-\la_1-\half)\Lambda_0^c + \sum_{k=1}^{j-1}
(\la_k-\la_{k+1})\Lambda_k^c + (\la_j-1)\Lambda_j^c +
\Lambda_{2l-j+1}^c,
$$
if $\la = (\la_1,\dots,\la_j,0,\dots,0)\otimes\det$.

According to \cite[Section~6.2]{W1}, there exist commuting actions
of $\cinf$ and of $O(2l+1)$ on $\F^{-l-\half}$.
\begin{prop} \cite[Theorem 6.2]{W1}
We have the following $(O(2l+1),\cinf)$-module decomposition:
$$
\F^{-l-\half} \cong \bigoplus_{\la\in\mathcal C_\half}
V_\la(O(2l+1))\otimes L(\cinf;\Lambda(\la),-l-\half),
$$
where $V_\la(O(2l+1))$ is the irreducible $O(2l+1)$-module
parameterized by $\la$.
%
\end{prop}

The operator $\Ct(t)$ acting on $\F^{-l-\half}$ can be expressed
now as
$$\begin{aligned}
\Ct(t) = &  \sum_{r\in\Z+\half}
 \sum_{p=1}^l t^r (-
\gamma_{-r}^{+,p}\gamma_r^{-,p} +
 \gamma_{-r}^{-,p}\gamma_r^{+,p})
 + \sum_{r\in\Z+\half} (-1)^{r+\hf}t^r\chi_{-r}\chi_r.
\end{aligned}
$$

\begin{definition}
The Bloch-Okounkov $n$-point $\cinf$-function of level $-l-\half$
(associated to a partition $\la$ of length $\leq l$) is defined as
$$
\mathfrak C_{-l-\half}^\la(q;t_1,\dots,t_n) = \trace_{L(\cinf;
\Lambda(\la),-l-\half)\oplus L(\cinf;
\Lambda(\la\otimes\det),-l-\half)} q^{L_0}\Ct(t_1)\cdots \Ct(t_n).
$$
\end{definition}

\begin{thm} \label{th:Clhalf}
The function $\mathfrak C^\la_{-l-\half}(q;t_1,\dots,t_n)$ is
given by
$$
\mathfrak C_{-\half}^{(0)} (q;t_1,\dots,t_n)\cdot\sum_{\sigma\in
W(B_l)} (-1)^{\ell(\sigma)}\mathfrak
C^{(k_1)}_{-1}(q;t_1,\dots,t_n)\cdots \mathfrak
C^{(k_l)}_{-1}(q;t_1,\dots,t_n),
$$
where $k_i = (\la+\rho_B-\sigma(\rho_B),\epsilon_i) \ge 0$.
\end{thm}

\begin{proof}
The proof is similar to that of Theorem \ref{cinf:thm:integral},
where instead we use the Weyl group of type $B_l$ and replace the
character $\mathsf{ch}_\la^o$ therein with
$$
\mathsf{ch}_\la^b(z_1,\dots,z_l) = \frac{\left|z_j^{\la_i +l -i
+\half} -z_j^{-(\la_i +l -i +\half)} \right|}{\left|z_j^{l-i
+\half} -z_j^{-(l -i +\half)}\right|}.
$$
Note that the factor $c_\la$ therein does not appear in this
computation, and also note that $z_i^{e_{ii}}$ does not act on
$\F^{-\half}$ which produces the appearance of the factor
$\mathfrak C^{(0)}_{-\half}(q;t_1,\dots,t_n)$ in the formula.
Finally $V_\la(O(2\ell+1))$ and $V_{\la\otimes{\rm
det}}(O(2\ell+1))$ are isomorphic as modules over the Lie algebra
of $O(2\ell+1)$, and hence they have the same character.
\end{proof}

The $1$-point $\cinf$-function of level $-l$ now follows by
combining Theorem~\ref{th:C1ptHalf}, Proposition~\ref{cinf:ctoa},
Theorem~\ref{th:Clhalf}, Remark~\ref{rem:residue}.

\subsection{The $q$-dimension of a $\cinf$-module of level $-l-\half$}

Again, let us denote by ${}^c\mathsf Q_{-l-\half}^\la(q)$ the
$q$-dimension of $L(\cinf;\Lambda(\la),-l-\half)\oplus
L(\cinf;\Lambda(\la\otimes\det),-l-\half)$.

\begin{thm}
The $q$-dimension ${}^c\mathsf Q_{-l-\half}^\la$ is
$$
{}^c\mathsf Q_{-l}^\la(q) =
\frac{1}{(q^\half)_\infty}\sum_{\sigma\in W(B_l)}
(-1)^{\ell(\sigma)} \;{}^c\mathsf Q_{-1}^{k_1}(q)\cdots {}^c\mathsf
Q_{-1}^{k_l}(q),
$$
where $k_i = (\la+\rho_B-\sigma(\rho_B),\varepsilon_i)$.
\end{thm}

\begin{proof}
This can be proved similarly to Theorem \ref{th:Clhalf} and we will
skip the details.
\end{proof}

\section{The $\dinf$-correlation functions and $q$-dimension formulas}
\label{sec:dinf}

As the methods involved in the $\dinf$ case are similar to the
cases of $\ainf$ and $\cinf$ treated in the previous sections, we
shall be more sketchy in this section.
\subsection{Lie algebra $\dinf$}
  \label{subsec_dinf}
Let
\begin{eqnarray*}
{\overline{d}}_{\infty} = \{ (a_{ij})_{i,j \in \Z} \in \gl \mid
a_{ij} = -a_{1-j,1-i} \}
 \end{eqnarray*}
be a Lie subalgebra of $\gl$ of type $D$ \cite{DJKM2}. Denote by
$\dinf = {\overline{d}}_{\infty} \bigoplus \C C $ the central
extension given by the restriction of $2$-cocycle (\ref{eq_cocy}) to
$\overline{d}_\infty$. Then $\dinf$ has a natural triangular
decomposition induced from $\hgl$ with Cartan subalgebra $(\dinf)_0
=  (\ainf)_{0} \cap \dinf$. Given $\Lambda \in (\dinf)_0^* $, we let
\begin{eqnarray*}
  H^d_i & = & E_{ii} + E_{-i, -i} - E_{i+1, i+1} - E_{-i+1, -i+1}
                                \quad (i \in \mathbb N),\\
  H^d_0 & = & E_{0,0} + E_{-1,-1} -E_{2,2} -E_{1,1} + 2C.
\end{eqnarray*}
Denote by $\hL^d_i $ the $i$-th fundamental weight of $\dinf$, i.e.
$\hL^d_i (H^d_j ) = \delta_{ij}$.  The Dynkin diagram of $\dinf$ is:
  
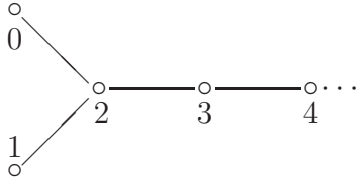
\begin{equation*}
  \begin{picture}(150,75) 
  \put(11,4){\line(1,1){25}}
  \put(5,60){$\circ$}
  \put(11,60){\line(1,-1){25}}
  \put(5,-1){$\circ$}
  \put(37,30){$\circ$}
  \put(44,33){\line(1,0){32}}
  \put(77,30){$\circ$}
  \put(84,33){\line(1,0){32}}
  \put(117,30){$\circ$}
  \put(124,30){$\cdots$}
  \put(5,48){$0$}
  \put(5,7){$1$}
  \put(38,20){$2$}
  \put(77,20){$3$}
  \put(117,20){$4$}
  \end{picture}
  \end{equation*}

\subsection{The $n$-point $\dinf$-functions of level $-l$}

Following \cite{W2, TW}, we introduce the following operator
$$
\D(t) = \sum_{k\in\N} (t^{k-\half}-t^{\half-k})(E_{k,k}-E_{1-k,1-k})
+ \frac{2}{t^\half-t^{-\half}}C.
$$
Regarding $\dinf$ as a subalgebra of $\ainf$, we have $\D(t) =
\A(t)-\A(t^{-1})$. According to \cite[Section 5.1]{W1} we have an
action of $\dinf$ on $\F^{-l}$, from which we obtain
$$
\D(t) = \sum_{p=1}^l\sum_{r\in\half+\Z} t^r(-
\gamma_{-r}^{+,p}\gamma_r^{-,p} +\gamma_{-r}^{-,p}\gamma_r^{+,p}).
$$
There is also an action of the Lie group $Sp(2l)$ on $\F^{-l}$ that
commutes with the action of $\dinf$.

Associated to a partition $\la =(\la_1,\ldots, \la_l)$ of length
$\leq l$, we define a highest weight $\Lambda(\la)$ for $\dinf$ to
be
$$
\Lambda(\la) = (-2l-\la_1-\la_2)\Lambda_0^d +
\sum_{k=1}^l(\la_k-\la_{k+1})\Lambda_k^d,
$$
with the convention that $\la_{l+1} = 0$.

\begin{prop} \cite[Theorem 5.2]{W1} \label{dualityCd}
We have the following $(\dinf,Sp(2l))$-module decomposition:
$$
\F^{-l} \cong \bigoplus_{\la} V_\la(Sp(2l))\otimes
L(\dinf;\Lambda(\la),-l),
$$
where $V_\la(Sp(2l))$ denotes the irreducible $Sp(2l)$-module of
highest weight $\la$.
\end{prop}

The $n$-point $\dinf$ function of level $-l$ (associated to a
partition $\la$ of length $\leq l$) is defined to be
$$
\mathfrak D^\la_{-l}(q;t_1,\dots,t_n) := \trace_{L(\dinf;
\Lambda(\la), -l)} q^{L_0} \D(t_1)\cdots \D(t_n),
$$
where $L_0$ is as usual a degree operator.

\begin{prop}
The function $\mathfrak D^{(m)}_{-1}(q;t_1,\dots,t_n)$ is given by
$$
\begin{aligned}
{[}z^m] & \sum_{\vec{\epsilon}_a\in\{\pm 1\}^n}
[\vec{\epsilon}_a]\trace_{\F^{-1}}z^{e_{11}}q^{L_0}\A(t_1^{\epsilon_1})\cdots
\A(t_n^{\epsilon_n})  \\
 & -[z^{m+2}]  \sum_{\vec{\epsilon}_a\in\{\pm 1\}^n}
[\vec{\epsilon}_a]\trace_{\F^{-1}}z^{e_{11}}q^{L_0}\A(t_1^{\epsilon_1})\cdots
\A(t_n^{\epsilon_n}),
\end{aligned}
$$
where as before we denote $\vec{\epsilon} =
(\epsilon_1,\dots,\epsilon_n)$ and $[\vec{\epsilon}] =
\epsilon_1\epsilon_2\cdots\epsilon_n$.
\end{prop}

\begin{proof}
Recalling that $\D(t) = \A(t)-\A(t^{-1})$ when acting on $\F^{-1}$,
we have
$$
\begin{aligned}
\trace_{\F^{-1}} & z^{e_{11}}q^{L_0} \D(t_1)\cdots \D(t_n) \\
 &= \prod_{j=1}^n \trace_{\F^{-1}}
z^{e_{11}}q^{L_0}\left(\A(t_j) -\A(t_j^{-1})\right) \\
& = \sum_{\vec{\eps}\in\{\pm 1\}^n} [\vec{\epsilon}]
\trace_{\F^{-1}} z^{e_{11}}q^{L_0} \A(t_1^{\eps_1})
\A(t_2^{\eps_2})\cdots \A(t_n^{\eps_n}) \\ &= \sum_{m\in\Z} z^m
\sum_{\vec{\eps}\in\{\pm 1\}^n} [\vec{\epsilon}]
\trace_{\F_{(m)}^{-1}} q^{L_0} \A(t_1^{\eps_1})
\A(t_2^{\eps_2})\cdots \A(t_n^{\eps_n}).
\end{aligned}
$$
Applying the trace of $z^{e_{11}}q^{L_0}\D(t_1)\cdots\D(t_n)$ to
both sides of Proposition~\ref{dualityCd} when $l=1$ yields
$$
\trace_{\F^{-1}}z^{e_{11}}q^{L_0}\D(t_1)\cdots\D(t_n) = \sum_{m \in
\Z_+} \mathsf{ch}_{(m)}^{sp}(z)\cdot
\trace_{L(d_\infty;\Lambda(m),-1)} q^{L_0}\D(t_1)\cdots\D(t_n).
$$
Combining the above equations, substituting the character formula
$\mathsf{ch}_{(m)}^{sp}(z) = (z^{m+1}-z^{-(m+1)})/(z-z^{-1})$, and
clearing denominators give us
$$\begin{aligned}
(z-z^{-1})& \sum_{m\in\Z} z^m \sum_{\vec{\eps}\in\{\pm 1\}^n}
[\vec{\epsilon_a}] \trace_{\F^{-1}}
q^{L_0} \A(t_1^{\eps_1}) \A(t_2^{\eps_2})\cdots \A(t_n^{\eps_n}) \\
&= \sum_{k\in\Z_{\geq 0}} (z^{k+1}-z^{-(k+1)})\cdot
\trace_{L(d_\infty;\Lambda(k),-1)} q^{L_0}\D(t_1)\cdots \D(t_n).
\end{aligned}
$$
Now the proposition follows.
\end{proof}

\begin{thm}
The function $\mathfrak D^\la_{-l}(q;t_1,\dots,t_n)$ is given by
$$
\sum_{\sigma\in W(C_l)} (-1)^{\ell(\sigma)}\mathfrak
D^{(k_1)}_{-1}(q;t_1,\dots,t_n)\cdots \mathfrak
D^{(k_l)}_{-1}(q;t_1,\dots,t_n),
$$
where $k_i = (\la+\rho_C-\sigma(\rho_C),\varepsilon_i)$.
\end{thm}

\begin{proof}
Recall the $Sp(2l)$-character formula (cf. \cite[24.18]{FH})
$$
\mathsf{ch}_\la^{sp}(z_1,\dots,z_l) = \frac{\left| z_j^{\la_i+l-i+1}
-z_j^{-(\la_i+l-i+1)}\right |}{\left |z_j^{l-i+1}-z_j^{-(l-i+1)}
\right |}.
$$
The proof proceeds in the same way as that of Theorem
\ref{gl:npt:lvll}, using now Proposition~\ref{dualityCd}, replacing
$\mathsf{ch}_\la^{gl}$ with $\mathsf{ch}_\la^{sp}$, clearing
denominators and comparing coefficients.  Note that $S_l$ therein is
replaced with the Weyl group $W(C_l)$.
\end{proof}
\subsection{The $q$-dimension of a $\dinf$-module of level $-l$}

Denote by ${}^d \mathsf Q_{-l}^\la(q)$ the $q$-dimension of
$L(\dinf;\Lambda(\la);-l)$.

\begin{prop}
Let $k \in \Z_+$. The $q$-dimension of the irreducible
$\dinf$-module of highest weight $\Lambda(k)$ and level $-1$ is
$$
{}^d\mathsf Q_{-1}^{(k)}(q) = \frac{1}{(q)_\infty^2}\sum_{m\geq 0}
(-1)^m q^{\hf m(m+1)} \left( q^{k(m+\half)}-
q^{(k+2)(m+\half)}\right).
$$
\end{prop}

\begin{proof}
This proof is similar to that of Proposition~\ref{gl:qdim:1}. By
Proposition~\ref{dualityCd} and the $Sp(2)$-character formula (and
clearing the Weyl denominator $(z-z^{-1})$), we arrive at the
following identity
$$
\sum_{r\in \Z_+} (z^{r+1}-z^{-(r+1)})\, {}^d\mathsf
Q^{(r)}_{-1}(q) = (z-z^{-1}) \frac1{(q)_{\infty}^2}  \sum_{r\in
\Z}
 \sum_{m=0}^\infty  z^r (-1)^m q^{\hf m(m+1)} q^{|r|(m+\hf)}.
 $$
The result now follows by comparing the coefficient of $z^{k+1}$
on both sides.
\end{proof}

\begin{thm}
The $q$-dimension of the irreducible $\dinf$-module of highest
weight $\Lambda(\la)$ and level $-l$ is
$$
{}^d\mathsf Q_{-l}^\la(q) = \sum_{\sigma\in W(C_l)}
(-1)^{\ell(\sigma)}\, {}^d\mathsf Q_{-1}^{(k_1)}(q)\cdots {}^d
\mathsf Q_{-1}^{(k_l)}(q),
$$
where $k_i = (\la+\rho_C-\sigma(\rho_C),\varepsilon_i)$.
\end{thm}

\begin{proof}
This follows from an appropriate change of the Weyl groups from
$S_l$ to $W(C_l)$ in the proof of Theorem \ref{gl:qdim:l}.
\end{proof}
\subsection{The $n$-point $\dinf$-functions of level $-l+\half$}

Introduce a neutral fermionic field $\phi (z) = \sum_{n \in \Z +\hf}
\phi_n z^{-n -\hf}$ which satisfies the following commutation
relations:
$$[ \phi_m , \phi_n ]_{+} = \delta_{m, -n},
     \quad m,n \in \Z +\hf.
$$
Denote by $\F^{\hf}$ the Fock space of $\phi (z)$. According to
\cite[Section~6.1]{W1}, the Fock space $\F^{-l+\hf} =\F^{-l}
\otimes \F^\hf$ admits the commuting actions of
$\mathfrak{osp}(1,2l)$ and $\dinf$.

Associated to a partition $\la =(\la_1,\ldots, \la_l)$ of length
$\leq l$, we define a highest weight $\Lambda(\la)$ for $\dinf$:
$$
\Lambda(\la) = (-2l+1-\la_1-\la_2)\Lambda_0^d +
\sum_{k=1}^l(\la_k-\la_{k+1})\Lambda_k^d,
$$
where $\la$ is a partition and we again take the convention that
$\la_{l+1} = 0$.
\begin{prop} \cite[Theorem 6.1]{W1}
We have the following $(\mathfrak{osp}(1,2l),\dinf)$-module
decomposition
$$
\F^{-l+\half} \cong \bigoplus_{\la}
V_\la(\mathfrak{osp}(1,2l))\otimes L(\dinf;\Lambda(\la),-l+\half),
$$
where the summation is over all partitions $\la$ with $\ell(\la)\le
l$.
\end{prop}

We can express the operator $\D(t)$ acting on $\F^{-l+\hf}$ as
$$
\D(t) = \sum_{r\in\half+\Z} \sum_{p=1}^l t^r  (-
\gamma_{-r}^{+,p}\gamma_r^{-,p} + \gamma_{-r}^{-,p}\gamma_r^{+,p})
+ \sum_{r\in\half+\Z} t^r\phi_{-r}\phi_r.
$$

Recall that the $n$-point $\dinf$-function of level $\hf$,
$\mathfrak D^{(0)}_\half(q;t_1,\dots,t_n)$, has been computed in
\cite{TW, W2}.
\begin{thm}
The $n$-point $\dinf$-function of level $-l+\hf$, $\mathfrak
D^\la_{-l+\half}(q;t_1,\dots,t_n)$, is given by
$$
\mathfrak D^{(0)}_{\half}(q;t_1,\dots,t_n) \sum_{\sigma\in W(B_l)}
(-1)^{\ell(\sigma)} \mathfrak D^{(k_1)}_{-1}(q;t_1,\dots,t_n)\cdots
\mathfrak D^{(k_l)}_{-1}(q;t_1,\dots,t_n),
$$
where $k_i = (\la+\rho_B-\sigma(\rho_B),\epsilon_i)$.
\end{thm}

\begin{proof}
The proof is similar to that of Theorem \ref{gl:npt:lvll}, where we
replace $\mathsf{ch}_\la^{gl}$ with $\mathsf{ch}_\la^{osp}$ (using
Lemma \ref{cinf:charequiv}).  Note that $z_i^{e_{ii}}$ does not act
on $\F^{\half}$ which produces the appearance of the factor
$\mathfrak D^{(0)}_{\half}(q;t_1,\dots,t_n)$ in the result.
\end{proof}
\subsection{The $q$-dimension of a $\dinf$-module of level $-l+\half$}

\begin{thm}
The $q$-dimension of the irreducible $\dinf$-module of highest
weight $\Lambda(\la)$ and level $-l+\half$ is
$$
{}^d\mathsf Q_{-l+\half}^\la(q) = (-q^\half)_\infty\sum_{\sigma\in
W(B_l)} (-1)^{\ell(\sigma)}\, {}^d\mathsf Q_{-1}^{(k_1)}(q)\cdots
{}^d \mathsf Q_{-1}^{(k_l)}(q),
$$
where $k_i = (\la+\rho_B-\sigma(\rho_B),\varepsilon_i)$.
\end{thm}

\begin{proof}

The proof proceeds as that of Theorem \ref{gl:qdim:l} with a few
changes. First, we substitute the Weyl group of type $B$ for $S_l$.
Now note that $z_i^{e_{ii}}$ does not act on $\F^{\half}$, which
produces a factor of $\trace_{\F^\half} q^{L_0}$ out front, which is
equal to $(-q^\half)_\infty$.
\end{proof}

\end{document}